\documentclass[draft]{agujournal2019}
\usepackage{url} %this package should fix any errors with URLs in refs.
\usepackage{lineno}
\usepackage[inline]{trackchanges}

\usepackage{soul}
\newcommand{\pd}[3]{\frac{\partial^{#3} #1}{\partial #2^{^{#3}}}}

\usepackage{xcolor}
\usepackage{amsmath}
\usepackage{amsfonts}
\usepackage{booktabs}
\usepackage{multirow}
\usepackage{pifont}

\draftfalse

\journalname{Journal of Advances in Modeling Earth Systems (JAMES)}

\begin{document}

%% ------------------------------------------------------------------------ %%
%  Title
%
% (A title should be specific, informative, and brief. Use
% abbreviations only if they are defined in the abstract. Titles that
% start with general keywords then specific terms are optimized in
% searches)
%
%% ------------------------------------------------------------------------ %%

% Example: \title{This is a test title}

\title{Evaluation of Implicit-Explicit Additive Runge-Kutta Integrators for the HOMME-NH Dynamical Core}

%% ------------------------------------------------------------------------ %%
%
%  AUTHORS AND AFFILIATIONS
%
%% ------------------------------------------------------------------------ %%

% Authors are individuals who have significantly contributed to the
% research and preparation of the article. Group authors are allowed, if
% each author in the group is separately identified in an appendix.)

% List authors by first name or initial followed by last name and
% separated by commas. Use \affil{} to number affiliations, and
% \thanks{} for author notes.
% Additional author notes should be indicated with \thanks{} (for
% example, for current addresses).

% Example: \authors{A. B. Author\affil{1}\thanks{Current address, Antartica}, B. C. Author\affil{2,3}, and D. E.
% Author\affil{3,4}\thanks{Also funded by Monsanto.}}

\authors{
Christopher J. Vogl\affil{1},
Andrew Steyer\affil{2},
Daniel R. Reynolds\affil{3},
Paul A. Ullrich\affil{4}, and
Carol S. Woodward\affil{1}
}

% \affiliation{1}{First Affiliation}
% \affiliation{2}{Second Affiliation}
% \affiliation{3}{Third Affiliation}
% \affiliation{4}{Fourth Affiliation}

\affiliation{1}{Center for Applied Scientific Computing, Lawrence Livermore National Laboratory, 7000 East Ave., Livermore, CA 94550, USA}
\affiliation{2}{Sandia National Laboratories, Albuquerque, NM  87185, PO Box 5800, MS 1320, USA}
\affiliation{3}{Department of Mathematics, Southern Methodist University, PO Box 750156, Dallas, TX 75275, USA}
\affiliation{4}{Department of Land, Air and Water Resources, University of California, Davis, Davis, CA 95616, USA}

%(repeat as many times as is necessary)

%% Corresponding Author:
% Corresponding author mailing address and e-mail address:

% (include name and email addresses of the corresponding author.  More
% than one corresponding author is allowed in this LaTeX file and for
% publication; but only one corresponding author is allowed in our
% editorial system.)

% Example: \correspondingauthor{First and Last Name}{email@address.edu}

\correspondingauthor{Christopher J. Vogl}{vogl2@llnl.gov}

%% Keypoints, final entry on title page.

%  List up to three key points (at least one is required)
%  Key Points summarize the main points and conclusions of the article
%  Each must be 100 characters or less with no special characters or punctuation

% Example:
% \begin{keypoints}
% \item	List up to three key points (at least one is required)
% \item	Key Points summarize the main points and conclusions of the article
% \item	Each must be 100 characters or less with no special characters or punctuation
% \end{keypoints}

\begin{keypoints}
\item When using certain ARK IMEX methods, the hydrostatic timestep size can be used for the non-hydrostatic model.
\item No single ARK IMEX method outperforms the rest in all metrics, but a small subset does excel in most metrics.
\item Hyperviscosity and vertical remap errors can affect global energy error with increased spatial resolution.
\end{keypoints}

%% ------------------------------------------------------------------------ %%
%
%  ABSTRACT
%
% A good abstract will begin with a short description of the problem
% being addressed, briefly describe the new data or analyses, then
% briefly states the main conclusion(s) and how they are supported and
% uncertainties.
%% ------------------------------------------------------------------------ %%

%% \begin{abstract} starts the second page

\begin{abstract}
  The nonhydrostatic High Order Method Modeling Environment (HOMME-NH) atmospheric dynamical core supports acoustic waves that propagate significantly faster than the advective wind speed, thus greatly limiting the timestep size that can be used with standard explicit time-integration methods.  Resolving acoustic waves is unnecessary for accurate climate and weather prediction.  This numerical stiffness is addressed herein by considering implicit-explicit additive Runge-Kutta (ARK IMEX) methods that can treat the acoustic waves in a stable manner without requiring implicit treatment of non-stiff modes.  Various ARK IMEX methods are evaluated for their efficiency in producing accurate solutions, ability to take large timestep sizes, and sensitivity to grid cell length ratio.  Both the Gravity Wave test and Baroclinic Instability test from the 2012 Dynamical Core Model Intercomparison Project (DCMIP) are used to recommend 5 of the 27 ARK IMEX methods tested for use in HOMME-NH.
\end{abstract}

\section*{Plain Language Summary}
  The Energy Exascale Earth System Model (E3SM) is an ongoing effort to produce actionable projections of variability and change in the Earth system.  Advances in computational power has allowed endeavors like E3SM to simulate physical phenomena at finer resolution than ever before.  Modeling the atmosphere at these finer resolutions requires improving both the underlying mathematical models and computational techniques.  This works focuses on the latter by evaluating cutting-edge computational techniques on a recently-developed atmosphere model and suggesting which of those techniques should be included in the next generation of E3SM.

%% ------------------------------------------------------------------------ %%
%
%  TEXT
%
%% ------------------------------------------------------------------------ %%

%%% Suggested section heads:
% \section{Introduction}
%
% The main text should start with an introduction. Except for short
% manuscripts (such as comments and replies), the text should be divided
% into sections, each with its own heading.

% Headings should be sentence fragments and do not begin with a
% lowercase letter or number. Examples of good headings are:

% \section{Materials and Methods}
% Here is text on Materials and Methods.
%
% \subsection{A descriptive heading about methods}
% More about Methods.
%
% \section{Data} (Or section title might be a descriptive heading about data)
%
% \section{Results} (Or section title might be a descriptive heading about the
% results)
%
% \section{Conclusions}

\section{Introduction}
  In recent years, the global atmospheric modeling community has directed substantial effort towards the development of nonhydrostatic global atmospheric modeling systems \cite{ullrich_dcmip2016:_2017}.  Although global weather prediction systems, such as the IFS \cite{wedi_modelling_2015}, have been run at nonhydrostatic resolutions for years, modern supercomputers are reaching the point where globally uniform nonhydrostatic simulations at sub-10km horizontal grid resolution are now possible on longer timescales for seasonal to decadal prediction \cite{reale_tropical_2017}.  Even on climatological time scales, investments in the development of regionally-refined models, also known as variable-resolution models, have permitted limited areas of the domain to be simulated with horizontal resolution that is more typical of regional climate models \cite{zarzycki_aquaplanet_2014, rauscher_impact_2014, harris_high-resolution_2016, huang_evaluation_2016, rhoades_sensitivity_2018}.

  A major hurdle in transitioning from hydrostatic to nonhydrostatic models is the numerical stiffness associated with vertical propagation of sound waves.  Such waves are largely irrelevant to weather or climate modeling, but the fine vertical grid spacing (on the order of tens of meters)  necessitated by the highly stratified nature of the atmosphere relative to the horizontal grid spacing (on the order of 1km or more) provides a major source of numerical stiffness.  The standard strategy of using fully implicit methods for integrating numerically stiff equations can be a computational burden in operational models \cite{evans_performance_2017, lott_algorithmically_2015} and consequently several alternatives have been developed.  One alternative to using fully implicit methods is to use modified equation sets that do not support vertically propagating waves \cite{ogura_scale_1962, durran_improving_1989, arakawa_unification_2009}.  However, these equations either cannot be employed on all scales or require global communication at each timestep (i.e. \citeA{davies_validity_2003, klein_regime_2010}).  A second and commonly used alternative among operational nonhydrostatic models \cite{ullrich_dcmip2016:_2017} is the use of implicit-explicit additive Runge-Kutta (ARK IMEX) methods \cite{ascher_implicit-explicit_1997, ullrich_operator-split_2012, weller_rungekutta_2013, gardner_implicitexplicit_2018}.  These methods work by distinguishing both ``slow'' and ``fast'' terms in the spatially discrete fluid equations, the latter set including vertically-propagating sound waves.  The ARK IMEX methods treat the ``fast'' terms implicitly while retaining an explicit treatment of the ``slow" terms, with coefficients and ordering devised in such a manner as to ensure stability and accuracy.  A vast library of ARK IMEX methods are now available throughout the published literature, but only a few studies have assessed these methods in the context of nonhydrostatic global atmospheric models \cite{gardner_implicitexplicit_2018, giraldo_implicit-explicit_2013}.  It is thus the goal of this paper to outline a number of metrics that can be used in performing this assessment and apply these metrics in the context of a new nonhydrostatic dynamical core.

  The dynamical core being assessed in this work is the spectral element nonhydrostatic High-Order Method Modeling Environment: HOMME-NH \cite{taylor_energy_2019}.  The hydrostatic version of HOMME-NH is presently used by the Energy Exascale Earth System Model (E3SM) atmospheric component model \cite{e3sm_project_energy_2018} and the Community Atmosphere Model Spectral Element (CAM-SE) dynamical core \cite{dennis_cam-se:_2012,rasch_overview_2018}.  The spectral element method used in HOMME-NH has many desirable properties for atmospheric modeling including parallel scalability, flexibility, and accuracy \cite{marras_review_2016}.  Of particular importance for this study, the eigenvalues of the spatial SE operator are purely imaginary when explicit diffusion is disabled \cite{ullrich_impact_2018}.  This characteristic suggests that a desirable property of ARK IMEX methods is that the stability region of the explicit method encompasses as large a region of the imaginary axis as possible \citeA <i.e.>[]{kinnmark_one_1984}.  Although the recommendations for ARK IMEX methods reached at the end of this paper is based on only a single model, the authors believe these conclusions are likely applicable to any problem where the spatial operators for the explicitly-treated problem components have purely imaginary eigenvalues.

  The dynamical core studied in \citeA{gardner_implicitexplicit_2018} is similar to HOMME-NH in that both use the spectral element method and, therefore, have a spatial operator with purely imaginary eigenvalues.  Furthermore, both dynamical cores partition the computational grid such that every vertical column is placed on only one computational node.  This allows for the implicit treatment of vertical terms without inter-node communication.  The two dynamical cores do, however, use different model formulations for the atmosphere.  Unlike the formulation studied in \citeA{gardner_implicitexplicit_2018}, there are only two forcing terms responsible for vertical acoustic wave propagation in the HOMME-NH formulation.  As such, this work considers a single splitting that treats only those two terms implicitly.  No additional investigation into treating horizontal terms implicitly is conducted because  \citeA{gardner_implicitexplicit_2018} found such splitting incurs too much communication cost to be beneficial.  Instead, the focus herein is on investigating 18 new ARK IMEX methods specifically developed for systems with purely imaginary eigenvalues, in addition to 9 methods from \citeA{gardner_implicitexplicit_2018} that performed well on the formulation considered there.  This work also extends the acceptable solution accuracy criterion in \citeA{gardner_implicitexplicit_2018} to a criterion based on surface pressure.

  The remainder of this paper begins with a brief review of the nonhydrostatic model used by HOMME-NH in Section \ref{SEC:nonhydrostatic}.  Section \ref{SEC:numerical} then discusses the spatial discretization and hyperviscosity treatment, which has implications on time integrator performance, as well as the 27 time integration schemes that will be evaluated.  Section \ref{SEC:evaluation} presents the various evaluation criteria and results including verification, energy conservation, and time-to-solution tests.  Recommendations on time integration schemes for both HOMME-NH and similar conservative systems, with purely imaginary spatial eigenvalues, are provided in Section \ref{SEC:recommendations}.  That section also recommends some improvements outside of the time integration scheme, which will likely affect time integration performance, before the paper concludes in Section \ref{SEC:conclusions}.

\section{The HOMME Nonhydrostatic Atmosphere Model} \label{SEC:nonhydrostatic}
  The formulation of the nonhydrostatic primitive equations used in HOMME-NH is a modified version of that developed by \citeA{laprise_euler_1992}:
  \begin{eqnarray}
    &\pd{\mathbf{u}}{t}{} + (\nabla_{\eta} \times \mathbf{u} + 2\boldsymbol{\Omega})\times \mathbf{u} + \frac{1}{2}\nabla_{\eta}(\mathbf{u} \cdot \mathbf{u}) + \frac{d\eta}{dt} \pd{\mathbf{u}}{\eta}{} + \frac{1}{\rho}\nabla_{\eta} p = 0, \nonumber \\
    &\pd{w}{t}{} + \mathbf{u} \cdot \nabla_{\eta} w + \frac{d\eta}{dt} \pd{w}{\eta}{} + g(1- \mu) = 0, \quad \mu = \pd{p}{\eta}{}/\pd{\pi}{\eta}{}, \nonumber \\
    &\pd{\phi}{t}{} + \mathbf{u} \cdot \nabla_{\eta} \phi + \frac{d\eta}{dt} \pd{\phi}{\eta}{} - gw  = 0, \label{EQ:thetanh} \\
    &\pd{\Theta}{t}{} + \nabla_{\eta} \cdot (\Theta \mathbf{u}) + \frac{\partial}{\partial \eta}(\Theta \frac{d\eta}{dt})  = 0, \quad \Theta = \pd{\pi}{\eta}{}\theta, \nonumber \\
    &\frac{\partial}{\partial t} (\pd{\pi}{\eta}{}) + \nabla_{\eta} \cdot (\pd{\pi}{\eta}{} \mathbf{u}) + \frac{\partial}{\partial \eta} (\pi \frac{d\eta}{dt}) = 0, \nonumber
  \end{eqnarray}
  where $g$ is the acceleration due to gravity, $\rho$ the density, $p$ the pressure, $w$ the velocity component in the direction of the spherical radius (vertical velocity), $\mathbf{u}$ the remaining velocity components (horizontal velocity), $\phi$ the geopotential, $\theta$ the potential virtual temperature, and $\frac{\partial \pi}{\partial \eta}$ the hydrostatic pressure gradient.  .  After noting that the hydrostatic pressure $\pi$ is the result of vertically integrating $g \rho$, the hybrid mass-based vertical coordinate $\eta$ is defined implicitly: $\pi = A(\eta)p_0 + B(\eta)p_s$, where $p_0$ is the top of atmosphere pressure, $p_s$ is the surface pressure, and $(A,B)$ are linear functions of $\eta$ such that $(A,B)=(1,0)$ when $\eta$ represents the top of the atmosphere and $(A,B)=(0,1)$ when $\eta$ represents the surface.  The horizontal gradient, divergence, and curl operators in this pressure-based coordinate system are denoted $\nabla_\eta$, $\nabla_\eta \cdot$, and $\nabla_\eta \times$, respectively.  Further details regarding the definition, derivation, and boundary conditions of the HOMME-NH primitive equations can be found in \citeA{taylor_energy_2019}.

\section{Numerical Methods} \label{SEC:numerical}
  To solve \eqref{EQ:thetanh}, the unknown fields $\mathbf{u}$, $w$, $\phi$, $\Theta$, $\pd{\pi}{\eta}{}$ are discretely represented on a cubed-sphere computational grid.  In the horizontal direction, the unknown fields and associated horizontal-derivative operators (i.e. $\nabla_{\eta} \cdot$, $\nabla_{\eta}$, and $\nabla_{\eta} \times$) in \eqref{EQ:thetanh} are represented using a fourth-order, mimetic, spectral-element method \cite{taylor_compatible_2010}.  These discretized operators are mimetic in that $\nabla_{\eta} \cdot$ and $\nabla_{\eta}$ are adjoints in a discrete sense.  In the vertical direction, the unknown fields are staggered with some quantities expressed at grid cell midpoints and others at grid cell interfaces.  The midpoint quantities are $\mathbf{u}$, $\Theta$, and$\pd{\pi}{\eta}{}$.  The interface quantities are $w$ and $\phi$.  For vertical-derivative operators (i.e. $\partial/\partial \eta$), the second-order accurate SB81 \cite{simmons_energy_1981} approach is used.  The SB81 discretization is also mimetic and supports a discrete product rule for the discrete spatial derivative operator $\partial/\partial \eta$.  Note that these mimetic properties result in a discrete system that both conserves energy and, when linearized about a steady-state solution, has eigenvalues along the imaginary axis.  It should also be noted that the computational grid can optionally be remapped in the vertical direction to improve numerical stability; however, the current use of the parabolic spline method \cite{zerroukat_parabolic_2006} does not conserve energy in a discrete sense.  Additional spatial discretization details can be found in \citeA{taylor_energy_2019}.

  \subsection{Integration in time}

    After the discretization of the spatial-differential operators, \eqref{EQ:thetanh} is now a system of ordinary differential equations: $d\mathbf{q}/dt = \mathbf{f}(\mathbf{q})$, where $\mathbf{q}(t)$ is the vector of all unknown grid quantities.  To address the numerical stiffness caused by the presence of acoustic waves, a typical approach is to separate out the components of $\mathbf{f}(\mathbf{q})$ responsible for vertical acoustic wave propagation.  For some nonhydrostatic models, this is not a trivial task \cite{gardner_implicitexplicit_2018}.  For \eqref{EQ:thetanh}, however, the terms to be separated are those that couple the vertical velocity $w$, the geopotential $\phi$, and the non-dimensional quantity $\mu$.  Note that if the system is hydrostatic, with no vertically propagating acoustic waves, then $\mu$ will be identically $1$ because $p = \pi$.  This observation motivates an additive splitting, $\mathbf{f}(\mathbf{q}) = \mathbf{f}^E(\mathbf{q}) + \mathbf{f}^I(\mathbf{q})$, given by
    \begin{eqnarray}
      &\frac{d\mathbf{u}_i}{dt} = \underbrace{-\big((\nabla_\eta \times \mathbf{u} + 2 \mathbf{\Omega}) \times \mathbf{u}\big)_i - \frac{1}{2} \big(\nabla_\eta(\mathbf{u}\cdot \mathbf{u})\big)_i - (\dot{\eta} \mathbf{u}_\eta)_i + (\nabla_\eta p / \rho)_i}_{\mathbf{f}^E_{\mathbf{u}_i}}, \nonumber \\
      &\frac{dw_j}{dt} = \underbrace{-(\mathbf{u} \cdot \nabla_\eta w)_j - (\dot{\eta}w_\eta)_j}_{f^E_{w_j}} + \underbrace{g(\mu_j-1)}_{f^I_{w_j}},
      \quad \frac{d\phi_j}{dt} = \underbrace{-(\mathbf{u} \cdot \nabla_\eta \phi)_j
      -  (\dot{\eta}\phi_\eta)_j}_{f^E_{\phi_j}} + \underbrace{gw_j}_{f^I_{\phi_j}}, \label{EQ:split_thetanh} \\
      &\frac{d \Theta_i}{dt} = \underbrace{-\big(\nabla_\eta \cdot (\Theta \mathbf{u})\big)_i - \big((\Theta \dot{\eta})_\eta\big)_i}_{f^E_{\Theta_i}},
      \quad \frac{d (\pi_\eta)_i}{dt} = \underbrace{-\big(\nabla_\eta \cdot (\pi_\eta \mathbf{u})\big)_i - \big((\pi_\eta \frac{d\eta}{dt})_\eta\big)_i}_{f^E_{(\pi_\eta)_i}}. \nonumber
    \end{eqnarray}
    where $i$ enumerates values at the center of the grid cell, vertically, and $j$ enumerates values at the vertical edges.

    ARK IMEX methods are designed for additively-split systems such as \eqref{EQ:split_thetanh}, where an explicit RK method is used on $\mathbf{f}^E$ while a diagonally-implicit RK method is used on $\mathbf{f}^I$.  With $\mathbf{q}^n$ denoting the approximate solution of $\mathbf{q}(t)$ at time $t_n$, these ARK IMEX methods approximate $\mathbf{q}(t_n+\Delta t)$, denoted $\mathbf{q}^{n+1}$, via
    \begin{align*}
      \mathbf{q}^{n+1} &= \mathbf{q}^n + \Delta t \sum_{i=1}^s\big(b^E_i \mathbf{f}^E(\mathbf{z}_i) + b^I_i \mathbf{f}^I(\mathbf{z}_i)\big) \text{, where} \\
      \mathbf{z}_i &= \mathbf{q}^n + \Delta t \sum_{j=1}^{i-1} a^E_{i,j} \mathbf{f}^E(\mathbf{z}_j) + \Delta t \sum_{j=1}^i a^I_{i,j} \mathbf{f}^I(\mathbf{z}_j), \quad i=1,\ldots,s.
    \end{align*}

    Various conditions on $a^E_{i,j}$, $a^I_{i,j}$, $b^E_i$, and $b^I_i$ exists to ensure a certain order of accuracy \cite{araujo_symplectic_1997}.  The remaining degrees of freedom can used to enforce stability properties and improve other aspects of the methods, leading to an abundance of potential IMEX RK methods in the literature.  Considering the purely-imaginary eigenvalues associated with the linearization of \eqref{EQ:thetanh} and the results of \citeA{gardner_implicitexplicit_2018}, 26 existing ARK IMEX methods and 1 explicit RK method are selected for evaluation.  In addition, a method developed by one of the authors specifically for nonhydrostatic atmosphere models is also evaluated.  All 28 methods are listed in Table \ref{TBL:methods}, along with their theoretical order of accuracy, number of stages requiring a nonlinear solve ($f^I$), number of stages requiring an evaluation of $\mathbf{f}^E(\mathbf{z}_j)$ ($f^E$), and references.  Table  \ref{tab:ark_properties} contains a more exhaustive list of properties for each ARK IMEX method, and plots of the stability regions for each method are available in a Zenodo archive \cite{vogl_stability_2019}.
    \begin{table}
      \centering
      \begin{tabular}{l c c c c}
        \hline
        Name & O & $f^I$ & $f^E$ & Ref. \\
        \hline
        KGU35* & 3 & 0 & 5 & {\scriptsize GU Eq. 56} \\
        \hline
        ARS222 & 2 & 2 & 3 & {\scriptsize ARS97 Sec. 2.6} \\
        \hline
        ARS232 & 2 & 2 & 3 & {\scriptsize ARS97 Sec. 2.5} \\
        \hline
        GSA222 & 2 & 2 & 3 & {\scriptsize BRS Eq. 5.2} \\
        \hline
        SSP2232 & 2 & 2 & 3 & {\scriptsize RMC Tbl. 2} \\
        \hline
        IMKG232(a,b) & 2 & 2 & 3 & {\scriptsize SVTG App.} \\
        \hline
        IMKG242(a,b) & 2 & 2 & 4 & {\scriptsize SVTG App.} \\
        \hline
        IMKG252(a,b) & 2 & 2 & 5 & {\scriptsize SVTG App.} \\
        \hline
        IMKG243a & 2 & 3 & 4 & {\scriptsize SVTG App.} \\
        \hline
        IMKG253(a,b) & 2 & 3 & 5 & {\scriptsize SVTG App.} \\
        \hline
        IMKG254(a,b,c) & 2 & 4 & 5 & {\scriptsize SVTG App.} \\
        \hline
      \end{tabular}
      \begin{tabular}{l c c c c}
        \hline
        Name & O & $f^I$ & $f^E$ & Ref. \\
        \hline
        ARS233 & 3 & 2 & 3 & {\scriptsize ARS Sec. 2.4} \\
        \hline
        SSP3333b & 3 & 2 & 3 & {\scriptsize CGG $\beta=\frac{2}{3}$} \\
        \hline
        SSP3333c & 3 & 2 & 3 & {\scriptsize CGG $\beta=\frac{\sqrt{3}+3}{6}$}\\
        \hline
        ARK324 & 3 & 3 & 4 & {\scriptsize KC03 ARK3(2)4} \\
        \hline
        ARS343 & 3 & 3 & 4 & {\scriptsize ARS Sec. 2.7} \\
        \hline
        IMKG343a & 3 & 3 & 4 & {\scriptsize SVTG App.} \\
        \hline
        ARS443 & 3 & 4 & 4 & {\scriptsize ARS Sec. 2.8} \\
        \hline
        DBM453 & 3 & 4 & 5 & {\scriptsize \ref{APP:dbm}} \\
        \hline
        ARK436 & 4 & 5 & 6 & {\scriptsize KC03 ARK4(3)6}\\
        \hline
        ARK437 & 4 & 6 & 7 & {\scriptsize KC19 ARK4(3)7} \\
        \hline
        ARK548 & 5 & 7 & 8 & {\scriptsize KC19 ARK5(7)8} \\
        \hline
      \end{tabular}
      \caption{List of numerical integration methods evaluated, listed with name, theoretical overall order (O), number of stages requiring a nonlinear solve ($f^I$), number of stages requiring an evaluation of $\mathbf{f}^E(\mathbf{z}_j)$ ($f^E$),  and references (ARS: \citeA{ascher_implicit-explicit_1997}, BRS: \citeA{boscarino_all_2018}, CGG: \citeA{conde_implicit_2017}, GU: \citeA{guerra_high-order_2016}, KC03: \citeA{kennedy_additive_2001}, KC19: \citeA{kennedy_higher-order_2019}, RMC: \citeA{rokhzadi_optimally_2018}, SVTG: \citeA{steyer_efficient_2019}, *\textit{fully explicit}).}
      \label{TBL:methods}
    \end{table}

  \subsection{ARK IMEX implementation using ARKode}

    Each stage $\mathbf{z}_i$ of an ARK IMEX method requires solving a nonlinear system $\mathbf{F}(\mathbf{z}_i) = \mathbf{0}$.  The formation and solving of this system is done using the ARKode package of SUNDIALS \cite{hindmarsh_sundials:_2005,sundials-web}.  While ARKode can adaptively adjust timestep sizes by default, this feature is disabled so the timestep size is fixed.  The nonlinear system is solved using a Newton approach, where approximations $\mathbf{z}_i^{(m)}$ to stages $\mathbf{z}_i$ for forming $\mathbf{q}^{n+1}$ are iterated from an initial guess $\mathbf{z}_i^{(0)}$, chosen to be $\mathbf{q}^n$:
    \begin{equation*}
      \mathbf{z}_i^{(m+1)} = \mathbf{z}_i^{(m)} + \boldsymbol{\delta}^{(m+1)}, \quad
      \bigg(\mathbf{I} - \Delta t\, a_{i,i}^I \nabla \mathbf{f}^I(\mathbf{z}_i^{(m)})\bigg)\,\boldsymbol{\delta}_i^{(m+1)} = \mathbf{b}(\mathbf{z}_i^{(m)}).
    \end{equation*}
    Because $\mathbf{f}^I$ only involves vertical derivatives in \eqref{EQ:split_thetanh}, the linear system for $\boldsymbol{\delta}_i^{(m+1)}$ is block-diagonal, with one block per vertical column of the computational grid.  The resulting linear solve requires no parallel communication because each block is local to single processor.  Furthermore, each block is tridiagonal , involving only $\phi_j$ and $w_j$, and is thus solved using LAPACK's tridiagonal solver routines DGTTRF and DGTTRS.  The stopping criterion for ARKode is that
    \begin{equation*}
      R_i^{m+1} ||\boldsymbol{\delta}_i^{(m+1)}||_\text{WRMS} < \epsilon, \quad R_i^{m+1} = \max \left(0.3R_i^m, \frac{||\boldsymbol{\delta}_i^{(m+1)}||_\text{WRMS}}{||\boldsymbol{\delta}_i^{(m)}||_\text{WRMS}}\right),
    \end{equation*}
    where $R_i^0 = 1$ and $\epsilon = 10^{-1}$.  For a solution vector $\mathbf{q}$ with $N$ components, the weighted root mean squared norm $||\cdot||_\text{WRMS}$ is defined as
    \begin{equation*}
       ||\mathbf{q}||_\text{WRMS} = \left( \frac{1}{N}\sum_{k=1}^N \left (\frac{q_k}{\epsilon_r |v_k| + \epsilon_{a,k}} \right)^2 \right)^{1/2},
    \end{equation*}
    where $\epsilon_r=10^{-6}$, $\epsilon_{a,k} = 10\epsilon_r$ when $q_k$ is either $\mathbf{u}_k$ or $w_k$, $\epsilon_{a,k} = 10^5\epsilon_r$ when $q_k$ is $\phi_k$, $\epsilon_{a,k}=10^6\epsilon_r$ when $q_k$ is $\Theta_k$, and $\epsilon_{a,k} = \epsilon_r$ when $q_k$ is $\frac{\partial \pi}{\partial \eta}_k$.  Here, $v_k$ corresponds to the $k$th component of the previous time-step solution.

    SUNDIALS 3.1.2 was the current version available when this work began and, as such, is the version used herein.  The ARKode framework in SUNDIALS 3.1.2 assumes a global nonlinear system and applies the aforementioned Newton method and stopping criterion to all grid cells together.  While this approach does not take full advantage of the column structure of the HOMME-NH nonlinear system, the global Newton approach used herein typically required 2 to 3 iterations to reach the stopping criterion using the values of $\epsilon_r$ and $\epsilon_{a,k}$ mentioned above.  The authors note that more recent versions of SUNDIALS allow for a user-defined nonlinear solver object, where the Newton method and stopping criterion could be applied to each grid column individually. Given the low numbers of iteratinos needed, we expect the nonlinear solver approach used here is sufficiently efficient to allow for accurate integrator assessments.

  \subsection{Hyperviscosity} \label{SEC:hyperviscosity}

    Although the spectral-element method used in HOMME-NH conserves energy in a discrete sense, it also produces persistent nonphysical waves \cite{ullrich_impact_2018}.  To address this in HOMME-NH,  hyperviscosity is applied by adding fourth-order, horizontal derivatives to the right-hand side of the ODE system: $d\mathbf{u}_i/dt = \mathbf{f^E_{\mathbf{u}_i}}(\mathbf{q}) - \nu (\Delta^2 \mathbf{u})_i$, where $\nu$ is the hyperviscosity magnitude.  This approach stabilizes the nonphysical waves at the cost of adding a tunable amount of energy dissipation to the system.  The hyperviscosity operator is applied, in a sub-stepped fashion if necessary, after \eqref{EQ:split_thetanh} is advanced in time.  More specifically, the ARK IMEX update is performed first and followed by the hyperviscosity update:
    \begin{enumerate}
      \item $\mathbf{q}^{n+1}_0 = \mathbf{q}^n + \Delta t \sum_{i=1}^s \big(b^E_i \mathbf{f}^E(\mathbf{z}_i) + b^I_i\mathbf{f}^I(\mathbf{z}_i)\big)$
      \item $\mathbf{u}^{n+1}_k = \mathbf{u}^{n+1}_{k-1} - \delta t\, \nu \Delta^2 \mathbf{u}^{n+1}_{k-1}$ for $k=1,\ldots,K$, where $\delta t = \Delta t / K$
      \item $\mathbf{u}^{n+1} = \mathbf{u}^{n+1}_K$
    \end{enumerate}
    This approach does limit formal convergence to first-order in time both because of the first-order Euler time integration and the first-order operator splitting; however, this limitation is usually only seen at timestep sizes significantly below operational $\Delta t$ values.  Furthermore, this approach only requires $K$ computations of $\Delta^2 \mathbf{u}$ versus computing $\Delta^2 \mathbf{z}_{\mathbf{u}}$ at each stage of the ARK IMEX method.

\section{Evaluation of ARK Methods} \label{SEC:evaluation}
  The ARK IMEX methods described in Section \ref{SEC:numerical} are evaluated on test problems defined at the 2012 Dynamical Core Model Intercomparison Project (DCMIP), summarized by \citeA{ullrich_dynamical_2012}.  The tests are all conducted on the Quartz machine at Lawrence Livermore National Laboratory.  Each Quartz node has a 2.1GHz Intel Xeon E5-2695 v4 CPU with 36 cores and 128GB of memory connected to an Omni-Path interconnect.   To build the HOMME standalone code, the Intel 16.0.3 compiler is used with MPICH 3.2,  NetCDF-C 4.5.0, NetCDF-Fortran 4.4.4, HDF5 1.10.3, and SUNDIALS 3.1.2.  A test problem simulating non-orographic gravity waves on a small planet is first used to verify the implementation of the methods.  A second test problem simulating a baroclinic instability on a full planet is used to assess energy conservation of the methods and then to evaluate the methods on accurate solution efficiency and the ability to take large timestep sizes.  Both test problems are run using $216$ MPI ranks across $6$ computational nodes.  The HOMME-NH code and namelist files are available in a Zenodo archive \cite{vogl_code_2019}.

  \subsection{Non-orographic Gravity Wave Test}
    For the non-orographic gravity wave test problem, a planet of radius $1/125$ that of Earth with no rotation or surface topography is used.  With a prescribed zonal wind speed, pressure and temperature fields can be chosen so that the system initially is in a steady state.  The initial condition for this test uses these fields with a perturbation added to the potential temperature to generate a gravity wave.  Full details of the setup are available from \citeA{ullrich_dynamical_2012}.

\begin{figure}
      \centering
      \includegraphics[width=12cm]{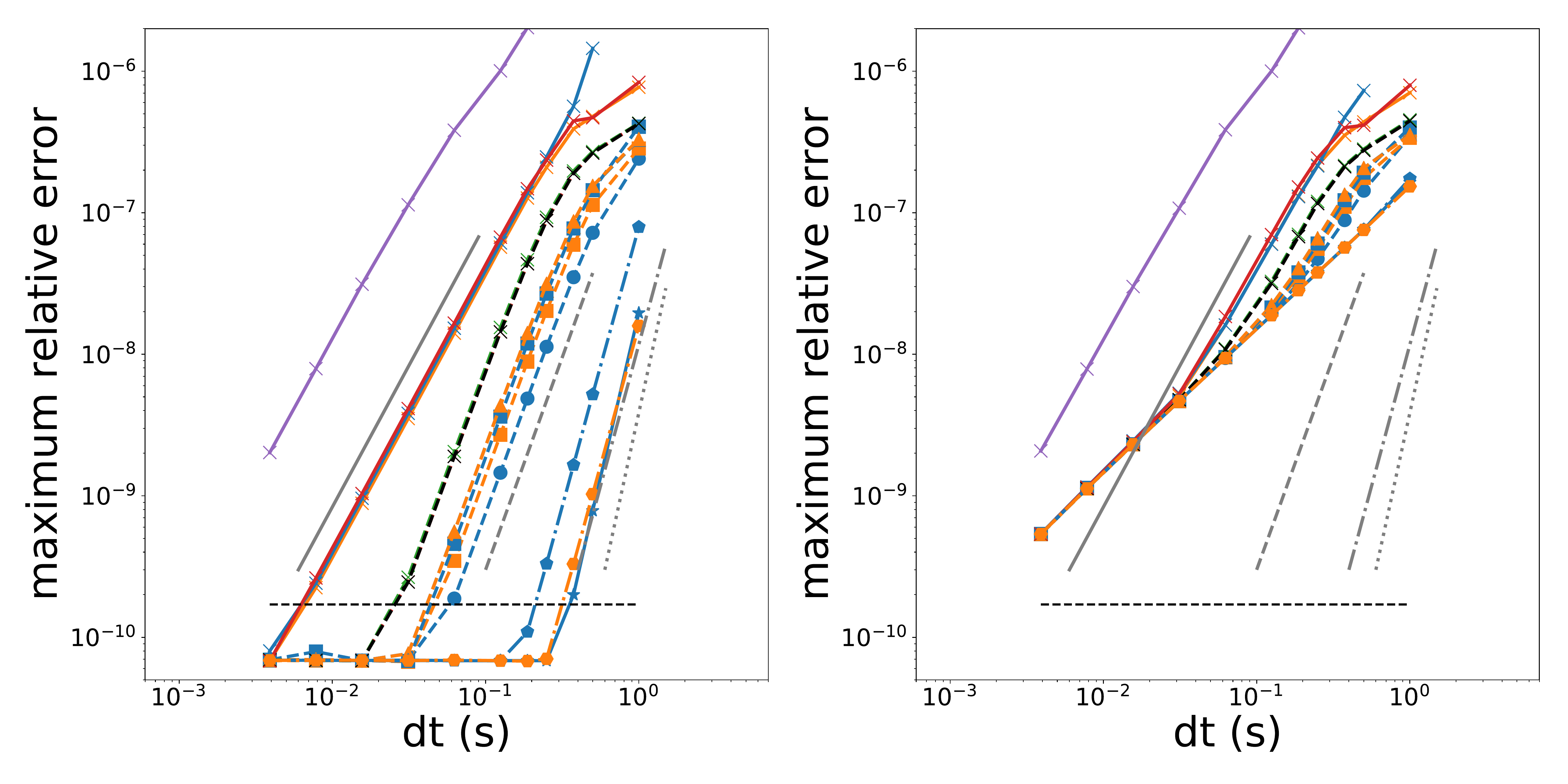}
      \caption{Convergence results for the gravity wave test showing the maximum relative error in temperature between various timestep sizes (s) and a reference solution (left: without hyperviscosity, right: with hyperviscosity, solid-lines: $2^\text{nd}$ or $5^\text{th}$ order methods, dash-lines: $3^\text{rd}$ order methods, dash-dot-lines: $4^\text{th}$ order methods; blue \textcolor{blue}{solid-line x}: ARS222, orange \textcolor{orange}{solid-line x}: ARS232, purple \textcolor{purple}{solid-line x}: GSA222, red \textcolor{red}{solid-line x}: SSP2232, red \textcolor{red}{dash-line x}: ARS233, green \textcolor{green}{dash-line x}: SSP3333b, black dash-line x: SSP3333c, orange \textcolor{orange}{dash-line triangle}: ARK324, blue \textcolor{blue} {dash-line triangle}: ARS343, blue \textcolor{blue}{dash-line square}: ARS443, orange \textcolor{orange}{dash-line square}: DBM453, blue \textcolor{blue}{dash-dot-line pentagon}: ARK436, orange \textcolor{orange}{dash-dot-line hexagon}: ARK437, blue \textcolor{blue}{solid-line star}: ARK548, right panel gray \textcolor{gray}{dotted-line}: first-order reference, left panel gray \textcolor{gray}{solid-line}: second-order reference, left panel gray \textcolor{gray}{dash-line}: third-order reference, left panel gray \textcolor{gray}{dash-dot-line}: fourth-order reference, left panel gray \textcolor{gray}{dotted-line}: fifth-order reference, black dash-dot-line: accumulated round-off error estimate).  A summary of convergence order without hyperviscosity is shown in Table \ref{TBL:gravitywave_convergence}}
      \label{FIG:gravitywave_convergence}
    \end{figure}

    A convergence study is performed to verify the implementation of the ARK IMEX methods, where observed convergence rates of the temperature field are compared to their theoretical predictions.  The spatial resolution is held fixed at $4,374$ fourth-order horizontal elements ($\approx1$km spacing) with $20$ vertical levels of uniform height for each element.  For this test, no vertical remap or hyperviscosity sub-stepping is used ($K=0$).  The timestep size ($\Delta t$) is varied for each of the ARK IMEX methods.  In lieu of an analytic solution, a reference solution from applying the explicit KGU35 method \cite{guerra_high-order_2016} with $\Delta t =  3.90625$E$-4$ is used.  The maximum relative temperature error across all grid cells at $t=5$h is shown for each $\Delta t$ in Figure \ref{FIG:gravitywave_convergence}.

    Note that in the absence of hyperviscosity, the errors of the various ARK IMEX methods decrease at a certain order as $\Delta t$ is decreased until they reach around $7*10^{-11}$.  The magnitude of this floor coincides with accumulated numerical round-off error.  Using the Intel \textit{epsilon} function, the accumulated round-off error is approximated by multiplying the obtained machine epsilon value ($2.220446049250313*10^{-16}$) by the number of timesteps taken in the reference solution ($768,000$).  The resulting value is shown in Figure \ref{FIG:gravitywave_convergence} as a horizontal dashed line.  A curve of the form $\alpha \Delta t^\beta$ is fit through the two lowest error values that reside above the dashed line for each method.  Table \ref{TBL:gravitywave_convergence} shows good agreement between the measured order $\beta$ and the order predicted by numerical analysis theory.  Similar convergence results for the IMKG methods can be found in \citeA{steyer_efficient_2019}.
    \begin{table}
      \centering
      \begin{tabular}{l c c}
        \hline
        Method & Order (T) & Order (M) \\
        \hline
        KGU35 & 3 & 2.95 \\
        \hline
        ARS222 & 2 & 1.99 \\
        \hline
        ARS232 & 2 & 1.98 \\
        \hline
        GSA222 & 2 & 1.97 \\
        \hline
        SSP2232 & 2 & 1.98 \\
        \hline
        ARS233 & 3 & 2.95 \\
        \hline
        SSP3333b & 3 & 2.95 \\
        \hline
        SSP3333c & 3 & 2.95 \\
        \hline
      \end{tabular}
      \hspace{1cm}
      \begin{tabular}{c c c}
        \hline
        Method & Order (T) & Order (M) \\
        \hline
        ARK324 & 3 & 2.96 \\
        \hline
        ARS343 & 3 & 2.96 \\
        \hline
        ARS443 & 3 & 2.99 \\
        \hline
        DBM453 & 3 & 2.96 \\
        \hline
        ARK436 & 4 & 3.97  \\
        \hline
        ARK437 & 4 & 3.95 \\
        \hline
        ARK538 & 5 & 4.75 \\
        \hline
      \end{tabular}
      \caption{Convergence test results for the gravity wave test in the absence of hyperviscosity (T: theoretical, M: measured)}
      \label{TBL:gravitywave_convergence}
    \end{table}

    Results with hyperviscosity applied are also shown in Figure \ref{FIG:gravitywave_convergence}.  The particular value for the hyperviscosity coefficient, $\nu = 5*10^8$ m$^4$/s, is chosen by tuning the value suggested by \citeA{ullrich_impact_2018} to the Gravity Wave test. As before, the errors of all the methods decrease as $\Delta t$ is decreased; however, they are now all limited to first order.  This is due to the time-split application of hyperviscosity, as discussed in Section \ref{SEC:numerical}.  Thus, the benefit of higher-order ARK IMEX methods at small $\Delta t$ is more-or-less lost, although higher-order methods do show an advantage at larger $\Delta t$.  The values of $\Delta t$ where higher-order ARK IMEX methods are advantageous will likely be problem-dependent; therefore, a test problem that incorporates more features of a full Earth system run is considered next.

  \subsection{Baroclinic Instability Test}
    Like the non-orographic gravity wave test, the baroclinic instability test prescribes a zonal wind field that is balanced by pressure and temperature fields.  A perturbation to the initial conditions is also added, only now it is to the zonal wind near the surface.  Unlike the linear evolution of the gravity wave test, the velocity field in this test becomes unsteady, resulting in a nonlinear evolution.  Full details of the setup are available from \citeA{ullrich_dynamical_2012}.  The spatial resolution used for all Baroclinic Instability tests is $5,400$ fourth-order horizontal elements ($\approx 110$km spacing) with $30$ vertical levels of varying heights chosen for more resolution near the surface.  For this test, hyperviscosity sub-stepping is used with $K=3$.

    \subsubsection{Global energy conservation}
      One metric for determining the ideal ARK IMEX method is how well the time integration scheme conserves global energy.  Recall from Section \ref{SEC:numerical} that while the mimetic discrete spatial operators conserve energy, the addition of hyperviscosity and vertical remapping procedure are not conservative and introduce energy conservation error.   Figure \ref{FIG:gravitywave_convergence} shows that the hyperviscosity error can dominate the overall error, making the time integration energy conservation error difficult to measure.  Furthermore, simulation of the baroclinic instability test for more than a few hours requires vertical remapping of the Lagrangian vertical coordinate.  Thus, short simulations are first conducted to compare the energy conservation error from each ARK IMEX method to the errors from the hyperviscosity and vertical remap approaches.  While one might use the initial condition of the baroclinic instability test itself for this, the evolution of those conditions is fairly linear for small time and not indicative of the nonlinear environment expected in HOMME-NH.  Thus, the result of simulating the baroclinic instability test for $7.5$ days with ARS232 at a timestep of $10$s, with hyperviscosity ($\nu = 10^{15}$ m$^4$/s) and vertical remap applied at each timestep, is used as an initial condition.

\begin{figure}[h]
        \centering
        \includegraphics[width=12cm]{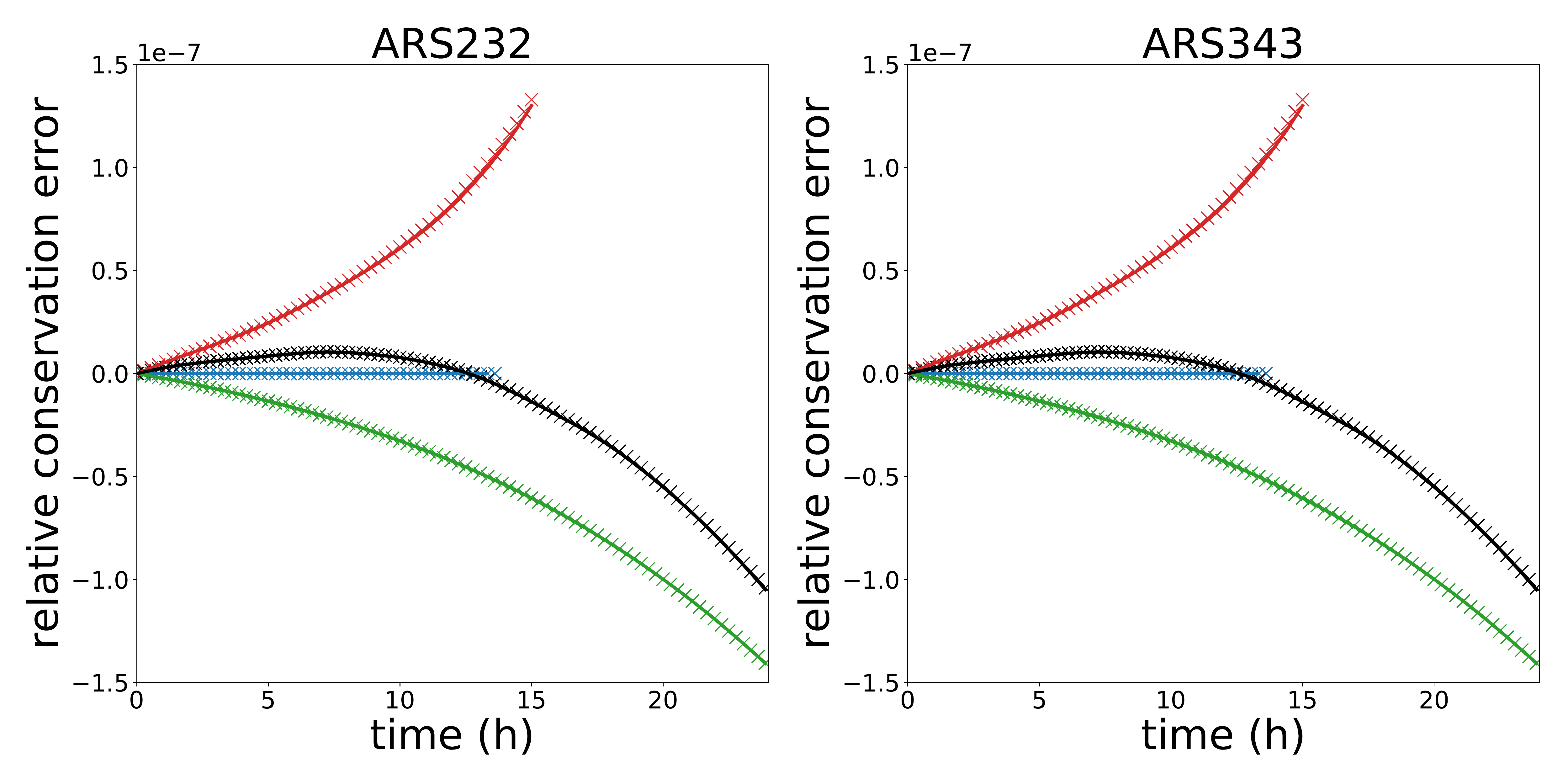}
         \caption{Measure of global energy conservation error from the numerical integrator, hyperviscosity, and vertical remap.  (left: ARS232, right: ARS343, \textcolor{blue}{blue}: integrator, \textcolor{green}{green}: integrator \& hyperviscosity, \textcolor{red}{red}: integrator \& vertical remap, black: integrator \& hyperviscosity \& remap, solid line: $\Delta t = 20$s, x-marker: $\Delta t =10$s)}
         \label{FIG:baroclinic_energy_error}
      \end{figure}

      Results are obtained for the second-order ARS232 and third-order ARS343 methods at timesteps of $10$s and $20$s.  These methods are run with four different setups: neither hyperviscosity nor vertical remap, with hyperviscosity ($\nu = 10^{15}$ m$^4/s$) but without vertical remap, without hyperviscosity but with vertical remap at each timestep, and with both hyperviscosity ($\nu = 10^{15}$ m$^4/s$) and vertical remap at each timestep.  Each setup is run for 24 hours or until the simulation goes unstable, usually indicated by a negative diagnosed density value.  Figure \ref{FIG:baroclinic_energy_error} shows the magnitude of relative error in the global energy $E(t)$ over time: $[E(t)-E(0)]/E(0)$.  Without hyperviscosity or vertical remap, the simulation goes unstable for both ARS232 and ARS343 in less than 14 hours.  The global energy conservation error for both methods at that point is at least an order of magnitude less than when hyperviscosity only is added or when vertical remap only is added.  It is worth noting that if either hyperviscosity or vertical remap are present, then neither using a higher-order method (ARS343 vs ARS232) or a smaller timestep ($10$s vs $20$s) results in an improvement in energy conservation error.  This lack of improvement indicates that the energy conservation errors from those sources is more dependent on the spatial resolution than on timestep size.  Another observation is that while it is expected that global energy is decreased by the addition of hyperviscosity, the addition of vertical remap without hyperviscosity increases global energy.  The implications of this behavior are discussed in Section \ref{SEC:recommendations}.

    \subsubsection{Largest accurate timestep size}
      Another metric for determining an ideal ARK IMEX method is how fast an accurate solution can be obtained.  This is, of course, dependent on the computational hardware, code compilation, definition of accuracy, and criterion for acceptability.  In addition to wall-clock run time, the largest accurate timestep size is analyzed in various ways to address the variability from computational hardware and code compilation.  The first measure is the largest accurate timestep size itself, which is independent of CPU clock speed or code compilation.  Another measure is the largest accurate timestep size normalized by the number of stages requiring a nonlinear solve, which addresses performance when solving the nonlinear system dominates the computational expense.  The final measure obtained is the largest accurate timestep normalized by the number of stages requiring an evaluation of $\mathbf{f}^E(\mathbf{z}_j)$, which addresses performance when inter-node communication dominates the computational expense.  Both hyperviscosity, with the $\nu$ parameter from Section \ref{SEC:hyperviscosity} set to $10^{15}$ m$^4$/s, and vertical remap are necessary to obtain a 15-day solution and are applied at every timestep.

\begin{figure}[h]
        \centering
        \includegraphics[width=12cm]{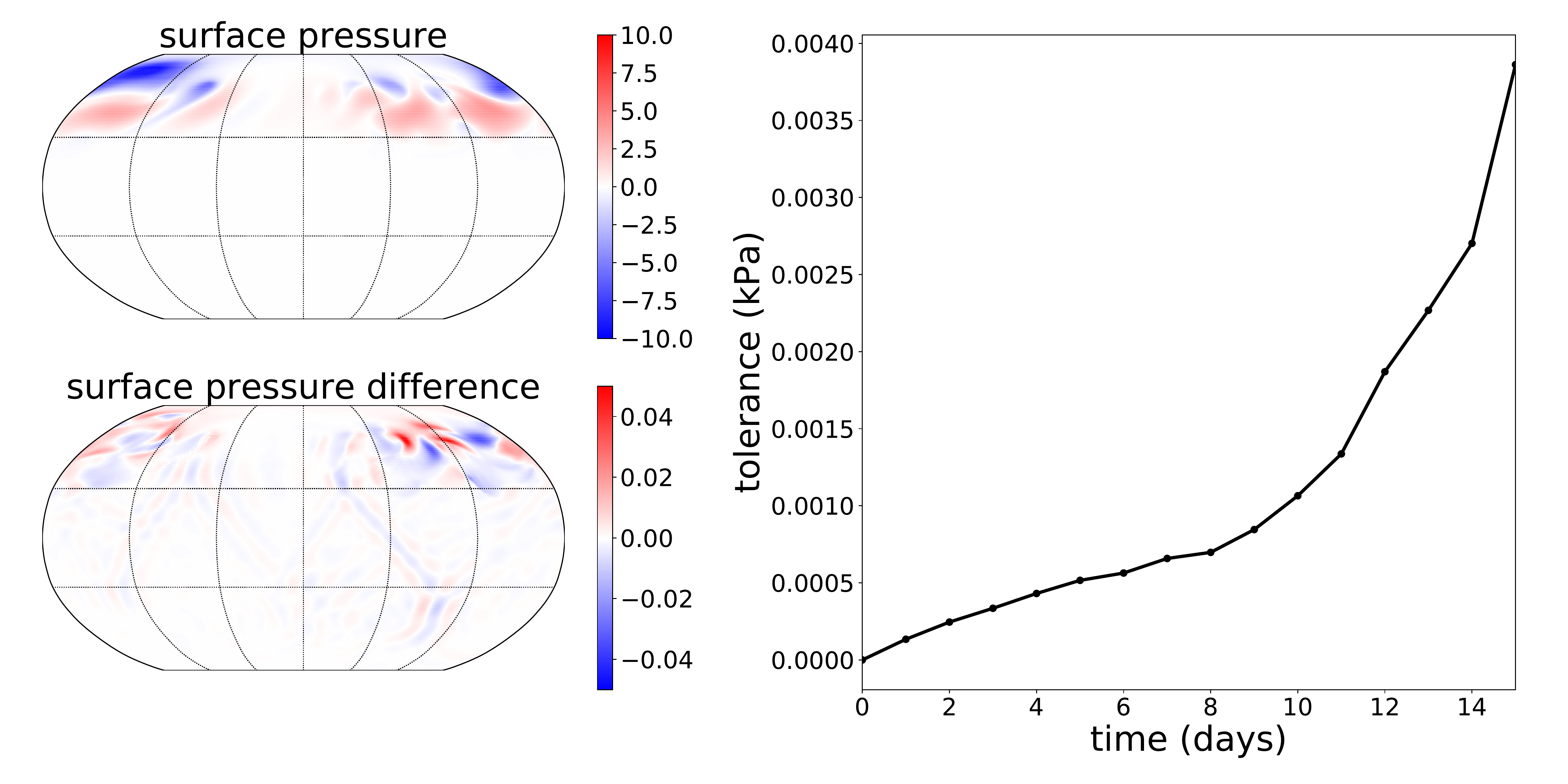}
        \caption{The surface pressure field (kPa) criterion used to determine largest accurate timestep size (top-left: reference surface pressure solution at 15 days, bottom-left: difference between reference solution at 15 days and initial condition, right: RMS difference over time that serves as a tolerance criterion).}
        \label{FIG:surface_pressure_tolerance}
      \end{figure}

      As in \citeA{jablonowski_baroclinic_2006} and \citeA{taylor_mass_2007}, the surface pressure field $p_s(\mathbf{s},t)$ is used to determine whether a solution is sufficiently accurate.  A reference surface pressure field is obtained at 24-hour intervals for 15 days by solving the hydrostatic version of \eqref{EQ:thetanh}, where $\mu$ is fixed to $1$, for the Baroclinic Instability test using KGU35 with a timestep of $10$s.  A second hydrostatic solution is obtained using KGU35 with the production timestep of $300$s.  The difference between the reference surface pressure field after 15 days and the initial condition used is shown in Figure \ref{FIG:surface_pressure_tolerance}.  Also shown is the difference between the reference and $\Delta t=300$s solutions after $15$ days and the RMS difference in surface pressure.  This RMS difference, denoted $\tau(t)$, is used as a tolerance for the accuracy of nonhydrostatic model \eqref{EQ:thetanh} solutions, which is valid because the Baroclinic Instability test is well within the hydrostatic regime when $5,400$ horizontal elements are used with $30$ levels.  Thus, the largest accurate timestep for each ARK IMEX method can now be defined as the maximum timestep of \{10, 20, 50, 100, 120, 135, 150, 160, 180, 192, 200, 216, 240, 270, 300, 320\} that results in a solution where the RMS difference of the surface pressure field from reference, denoted $\delta(t)$, is less than the tolerance $\tau(t)$ depicted in Figure \ref{FIG:surface_pressure_tolerance} for all $15$ days.

      The largest accurate timestep size and corresponding wall-clock run times are shown in Table \ref{TBL:baroclinic_accuracy_results}.  For the largest accurate timestep overall, the ARK548 method produces an accurate solution with the hydrostatic timestep $\Delta t=300$s.  For the largest accurate timestep normalized by number of stages requiring a nonlinear solve ($\Delta t / f^I$), IMKG242b (2 implicit stages) out performs the ARK548 method and produces an accurate solution in the least amount of time overall.  This is consistent with the environment in which these solutions are produced (216 ranks across 6 nodes), where obtaining the nonlinear solutions dominates the computational cost.  If enough computational nodes are used such that communication becomes the dominant cost, one might expect IMKG243a to reach an accurate solution before the other methods, as it has the largest accurate timestep normalized by number of stages requiring an evaluation of $\mathbf{f}^E(\mathbf{z}_j)$ ($\Delta t / f^E$).
      \begin{table}
        \centering
        \begin{tabular}{lcccc}
          \hline
          Method & $\Delta t$ & $\frac{\Delta t}{f^I}$ & $\frac{\Delta t}{f^E}$ & Time \\
          \hline
          ARS222 & 20 & 10 & 7 & 8.75 \\
          \hline
          ARS232 & 120 & 60 & 40 & 1.54 \\
          \hline
          GSA222 & 20 & 10 & 7 & 8.81 \\
          \hline
          SSP2232 & 135 & 68 & 45 & 1.49 \\
          \hline
          IMKG232a & 50 & 25 & 17 & 3.92 \\
          \hline
          IMKG232b & 50 & 25 & 17 & 3.93 \\
          \hline
          IMKG242a & 160 & \textbf{80} & 40 & 1.36 \\
          \hline
          IMKG242b & \textbf{240} & \textbf{120} & \textbf{60} & \textbf{0.90} \\
          \hline
          IMKG243a & \textbf{270} & \textbf{90} & \textbf{68} & \textbf{0.92} \\
          \hline
          IMKG252a & 100 & 50 & 20 & 2.38 \\
          \hline
          IMKG252b & 120 & 60 & 24 & 1.96 \\
          \hline
          IMKG253a & 100 & 33 & 20 & 2.72 \\
          \hline
          IMKG253b & 120 & 40 & 24 & 2.21 \\
          \hline
        \end{tabular}
        \hspace{1cm}
        \begin{tabular}{lcccc}
          \hline
          Method & $\Delta t$  & $\frac{\Delta t}{f^I}$ & $\frac{\Delta t}{f^E}$ & Time \\
          \hline
          IMKG254a & 120 & 30 & 24 & 2.43 \\
          \hline
          IMKG254b & 120 & 30 & 24 & 2.51 \\
          \hline
          IMKG254c & 100 & 25 & 20 & 2.99 \\
          \hline
          ARS233 & 180 & \textbf{90} & \textbf{60} & 1.07 \\
          \hline
          SSP3333$\frac{b}{c}$ & 180 & \textbf{90} & \textbf{60} & $\frac{1.07}{1.12}$ \\
          \hline
          ARK324 & \textbf{240} & \textbf{80} & \textbf{60} & \textbf{1.01} \\
          \hline
          ARS343 & 200 & 67 & 50 & 1.19 \\
          \hline
          IMKG343a & 160 & 53 & 40 & 1.54 \\
          \hline
          ARS443 & 135 & 34 & 34 & 2.09 \\
          \hline
          DBM453 & \textbf{270} & 68 & \textbf{54} & 1.09 \\
          \hline
          ARK436 & 216 & 43 & 36 & 1.61 \\
          \hline
          ARK437 & 50 & 8 &  7 & 7.59 \\
          \hline
          ARK548 & \textbf{300} & 50 &  43 & 1.5 \\
          \hline
        \end{tabular}
        \caption{Largest timestep size ($\Delta t$), in seconds, for each method that results in a surface pressure error below tolerance and corresponding run time, in hours ($f^I$: number of stages requiring a nonlinear solve, $f^E$: number of stages requiring an evaluation of $\mathbf{f}^E(\mathbf{z}_j)$).}
        \label{TBL:baroclinic_accuracy_results}
      \end{table}

    \subsubsection{Ability to use hydrostatic timestep size}
      A final metric for determining an ideal ARK IMEX method is whether the current hydrostatic dynamical core timestep of $\Delta t = 300$s can be used.  Certain components of an Earth system model might require HOMME-NH to advance the solution by the hydrostatic timestep before the solution is coupled back to those components.  As such, the ability to take the hydrostatic timestep is assessed using the setup of the previous section but without an accuracy criterion.  Here, ``exceedance'' is defined as how much the RMS deviation of surface pressure from the reference solution ($\delta(t)$) exceeds the RMS tolerance ($\tau(t)$). Table \ref{TBL:baroclinic_coupling_results} shows the largest timestep sizes that yield a 15 day solution, corresponding wall-clock run times, and relative maximum exceedance values, the latter of which is defined as
      \begin{gather*}
        \frac{\max_i \big(\delta(t_i)-\tau(t_i)\big)}{\delta(t_j)}, \quad \text{where} \quad j = \text{argmax}_j \big(\delta(t_i)-\tau(t_i)\big).
      \end{gather*}
      The six ARK IMEX methods that can produce a solution using $\Delta t=300$s are ARK548, DBM453, IMKG252b, IMKG253b, IMKG254a, and IMKG254b methods.  Of those six, IMKG252b has the largest timestep normalized by number of stages requiring a nonlinear solve and obtains the solution in the least amount of time.  As before, the method with the largest timestep normalized in this fashion is expected to excel in the environment where the nonlinear solution dominates computation cost.  Looking at the largest timestep normalized by number of stages requiring an evaluation of $\mathbf{f}^E(\mathbf{z}_j)$, one might expect equal performance across the DBM453, IMKG252b, IMKG253b, IMKG254a, and IMKG254b methods.  DBM453 shows significantly higher accuracy than the IMKG methods using $\Delta t =300$s while requiring similar run time as IMKG254a or IMKG254b and about 25\% more run time than IMKG252b.
      \begin{table}
        \centering
        \begin{tabular}{lccccc}
          \hline
          Method & $\Delta t$ & $\frac{\Delta t}{f^I}$ & $\frac{\Delta t}{f^E}$ & Time & Exc.\\
          \hline
          ARS222 & 50 & 25 & 17 & 3.57 & .24 \\
          \hline
          ARS232 & 120 & 60 & 40 & 1.54 & 0 \\
          \hline
          GSA222 & 50 & 25 & 17 & 3.62 & .25 \\
          \hline
          SSP2232 & 135 & 68 & 45 & 1.49 & 0 \\
          \hline
          IMKG232a & 100 & 50 & 33 & 1.96 & .57 \\
          \hline
          IMKG232b & 200 & 100 & 67 & 0.99 & .71 \\
          \hline
          IMKG242a & 160 & 80 & 40 & 1.36 & 0 \\
          \hline
          IMKG242b & 240 & 120 & 60 & 0.90 & 0 \\
          \hline
          IMKG243a & 270 & 90 & 68 & 0.92 & 0 \\
          \hline
          IMKG252a & 160 & 80 & 32 & 1.49 & .38 \\
          \hline
          IMKG252b & \textbf{300} & 150 & 60 & 0.79 & \textbf{.61} \\
          \hline
          IMKG253a & 200 & 67 & 40 & 1.38 & .55 \\
          \hline
          IMKG253b & \textbf{300} & 100 & 60 & 0.92 & \textbf{.61} \\
          \hline
        \end{tabular}
        \begin{tabular}{lccccc}
          \hline
          Method & $\Delta t$ & $\frac{\Delta t}{f^I}$ & $\frac{\Delta t}{f^E}$ & Time & Exc.\\
          \hline
          IMKG254a  & \textbf{300} & 75 & 60 & 0.98 & \textbf{.61} \\
          \hline
          IMKG254b  & \textbf{300} & 75 & 60 & 1.01 & \textbf{.61} \\
          \hline
          IMKG254c & 160 & 40 & 32 & 1.84 & .38 \\
          \hline
          ARS233 & 200 & 100 & 67 & 0.98 & .11 \\
          \hline
          SSP3333$\frac{b}{c}$ & 200 & 100 & 67 & $\frac{0.98}{0.96}$ & .11 \\
          \hline
          ARK324 & 240 & 80 & 58 & 1.01 & 0 \\
          \hline
          ARS343 & 200 & 67 & 50 & 1.19 & 0 \\
          \hline
          IMKG343a & 270 & 90 & 68 & 0.92 & .55 \\
          \hline
          ARS443 & 135 & 34 & 34 & 2.09 & 0 \\
          \hline
          DBM453 & \textbf{300} & 75 & 60 & 0.99 & \textbf{.03} \\
          \hline
          ARK436 & 216 & 43 & 36 & 1.61 & 0 \\
          \hline
          ARK437 & 50 & 8 & 7 & 7.59 & 0 \\
          \hline
          ARK548 & \textbf{300} & 50 &  43 & 1.5 & \textbf{0} \\
          \hline
        \end{tabular}
        \caption{Largest timestep size ($\Delta t$), in seconds, for each method that results in a solution after 15 days, the corresponding wall-clock run time, in hours, and the relative maximum exceedance of RMS tolerance ($f^I$: number of stages requiring nonlinear solves, $f^E$: number of stages requiring evaluation of $\mathbf{f}^E(\mathbf{z}_j)$).}
        \label{TBL:baroclinic_coupling_results}
      \end{table}

      Recall from Section \ref{SEC:numerical} that terms related to vertically-propagating acoustic waves are treated with the implicit Butcher Tableau while the remaining terms, including those related to horizontally-propagating acoustic waves are treated with the explicit Butcher Tableau.  If the horizontal grid resolution continues to improve while the vertical grid resolution remains constant, the stability properties of the explicit Butcher Tableau will become more important.  Thus, each ARK IMEX method is now used to obtain solutions on planets of decreasing radii to investigate the effect of shrinking horizontal grid cell lengths.  Two planets are considered: one with a radius $1/10$ that of Earth and a rotation rate $10$ times that of Earth, where the hydrostatic timestep size is $30$s, and one with a radius $1/100$ that of Earth and a rotation rate $100$ times that of Earth, where the hydrostatic timestep size is $3$s.  Following \citeA{ullrich_dynamical_2012}, these two planets are denoted ``X10'' and ``X100'' respectively.  With the number of grid cells held constant, these setups result in smaller relative horizontal grid cell lengths.  Note that the hyperviscosity is adjusted to reflect the different grid cell size ($\nu_\text{X10}=10^{12}$ m$^4$/s and $\nu_\text{X100}=10^9$ m$^4$/s).  If the stability of the ARK IMEX method is not sensitive to this change, the largest timestep size should scale with the radius.  The results in Table \ref{TBL:baroclinic_coupling_small_planet_results} show varying levels of stability sensitivity.  Of the methods that can produce a 15-day solution with the hydrostatic timestep ($\Delta t = 300$s), ARK548 shows strong sensitivity, requiring a timestep size of $19.2$s on the X10 planet and of $1.6$s on the X100 planet.  DBM453 shows more muted sensitivity, requiring a timestep size of $27$s on the X10 planet and a timestep size of $2.16$ on the X100 planet.  That said, the DBM453 timestep sizes are within 10\% of the $30$s hydrostatic timestep size that IMKG254a, IMKG254b, and IMKG252b can take on the X10 planet, within 2\% of the $2.2$s timestep size that IMKG252b and IMKG254a can take on the X100 planet, and within 10\% of the $2.4$s timestep size that IMKG254b can take on the X100 planet.
      \begin{table}
        \centering
        \begin{tabular}{lccc}
          \hline
          Method & $\Delta t$ & $\Delta t_\text{X10}$ & $\Delta t_\text{X100}$ \\
          \hline
          ARS222 & 50 & 5.0 & 0.5 \\
          \hline
          ARS232 & 120 & 12.0 & 1.20 \\
          \hline
          GSA222 & 50 & 5.0 & 0.50 \\
          \hline
          SSP2232 & 135 & 13.5 & 1.35 \\
          \hline
          IMKG232a & 100 & 10.0 & 1.00 \\
          \hline
          IMKG232b & 200 & 19.2 & 1.92 \\
          \hline
          IMKG242a & 160 & 16.0 & 1.60 \\
          \hline
          IMKG242b & 240 & 24.0 & 2.40 \\
          \hline
          IMKG243a & 270 & 27.0 & 2.16 \\
          \hline
          IMKG252a & 160 & 16.0 & 1.60 \\
          \hline
          IMKG252b & \textbf{300} & \textbf{30.0} & \textbf{2.20} \\
          \hline
          IMKG253a & 200 & 20.0 & 2.00 \\
          \hline
          IMKG253b & \textbf{300} & \textbf{30.0} & \textbf{2.16} \\
          \hline
        \end{tabular}
        \hspace{1cm}
        \begin{tabular}{lccc}
         \hline
          Method & $\Delta t$ & $\Delta t_\text{X10}$ & $\Delta t_\text{X100}$\\
          \hline
          IMKG254a & \textbf{300} & \textbf{30.0} & \textbf{2.16} \\
          \hline
          IMKG254b & \textbf{300} & \textbf{30.0} & \textbf{2.40} \\
          \hline
          IMKG254c & 160 & 16.0 & 1.60 \\
          \hline
          ARS233 & 200 & 18.0 & 1.80 \\
          \hline
          SSP3333$\frac{b}{c}$ & 200 & 18.0 & 1.80 \\
          \hline
          ARK324 & 240 & 19.2 & 1.60 \\
          \hline
          ARS343 & 200 & 19.2 & 1.60 \\
          \hline
          IMKG343a & 270 & 24.0 & 2.16 \\
          \hline
          ARS443 & 135 & 13.5 & 1.35 \\
          \hline
          DBM453 & \textbf{300} & \textbf{27.0} & \textbf{2.16} \\
          \hline
          ARK436 & 216 & 20.0 & 1.80 \\
          \hline
          ARK437 & 50 & 5.0 & 0.50 \\
          \hline
          ARK548 & 300 & 19.2 & 1.60 \\
          \hline
        \end{tabular}
        \caption{Largest timestep size, in seconds, for each method that results in a solution after 15 revolutions on planets of various radii ($\Delta t$: radius equal to that of Earth, $\Delta t_\text{X10}$: radius $1/10^\text{th}$ that of Earth, $\Delta t_\text{X100}$: radius $1/100^\text{th}$ that of Earth).}
        \label{TBL:baroclinic_coupling_small_planet_results}
      \end{table}

\section{Recommendations} \label{SEC:recommendations}
  In Section \ref{SEC:evaluation}, the 27 ARK IMEX methods listed in Section \ref{SEC:numerical} were evaluated on the metrics of largest timestep size that yields an acceptable accurate solution, largest timestep size that yields a 15-day solution of any accuracy, and sensitivity to horizontal grid resolution.  The error contributions of the hyperviscosity and vertical remap approaches of HOMME-NH were also investigated.  Noting that the HOMME-NH dynamical core will be solving fully nonlinear, stiff problems from earth system models, the Baroclinic Instability test from 2012 is primarily chosen to evaluate the ARK IMEX methods, as it is both stiff and exhibits nonlinear phenomena.  While the performance of any given ARK IMEX method will depend on the particular problem being solved, as well as the computational hardware, the authors are confident in making general recommendations for future development of both HOMME-NH and other dynamical cores.

  \subsection{ARK IMEX methods for HOMME-NH}
    The performance of the various ARK IMEX methods on the Baroclinic Instability test for each of the metrics is summarized below by metric, with brackets grouping ARK IMEX methods of equivalent performance within each metric:
    \begin{enumerate}
      \item Largest accurate timestep size ($\Delta t$): ARK548, \{DBM453, IMKG243a\}, \{IMKG242b, ARK324\}
      \item Efficiency in nonlinear solution cost dominated environments ($\Delta t/f^I$): IMKG242b, \{IMKG243a, ARS233, SSP3333(b,c)\}, and \{ARK324, IMKG242a\}
      \item Efficiency in communication cost dominated environments ($\Delta t/f^E$): IMKG243a, \{ARK324, ARS233, IMKG242b, SSP3333(b,c)\}, and DBM453
      \item Able to take hydrostatic timestep: \{ARK548, DBM453, IMKG252b, IMKG253b, IMKG254a, IMKG254b\}
      \item Most accurate solution with hydrostatic timestep: ARK548
      \item Able to take hydrostatic timestep on X10 planet: \{IMKG252b, IMKG253b, IMKG254a, IMKG254b\}
    \end{enumerate}
    While no one ARK IMEX method outperforms the rest in all metrics, there are a handful of methods that consistently perform better than the rest.  Thus, the authors recommend that the methods in Table \ref{TBL:recommendations}, listed in alphabetical order, be implemented and considered for operation in HOMME-NH, followed by accuracy assessments through more in-depth code verification tests.  ARK548 does not appear in Table \ref{TBL:recommendations} because while the method was able to produce an accurate solution at the hydrostatic timestep, the large number of both explicit and implicit stages make the cost of each timestep significantly higher than that of the methods in Table \ref{TBL:recommendations}.
    \begin{table}
      \centering
      \begin{tabular}{lccccc}
        & Accuracy & Solver Cost & Communication Cost & Stability & X10 Stability \\
        \hline
        ARK324 & \checkmark & \checkmark & \checkmark & & \\
        \hline
        DBM453 & \checkmark* & & \checkmark & \checkmark* &  \\
        \hline
        IMKG242b & \checkmark & \checkmark* & \checkmark & & \\
        \hline
        IMKG243a & \checkmark* & \checkmark & \checkmark* & & \\
        \hline
        IMKG252b & & & & \checkmark* & \checkmark* \\
        \hline
      \end{tabular}
      \caption{Recommended ARK IMEX methods along with performance in various metrics: largest accurate timestep size (accuracy), efficiency in nonlinear solution cost dominated environments (solver cost), efficiency in communication cost dominated environments (communication cost), ability to take hydrostatic timestep (stability), and ability to take hydrostatic timestep on X10 planet (X10 stability) (\checkmark: method outperformed most other methods, \checkmark*: method performed best of all methods)}
      \label{TBL:recommendations}
    \end{table}

  \subsection{Hyperviscosity and Vertical Remap}
    It is a major finding of this paper that the benefits of these time integration methods, namely high order-of-accuracy and energy conservation, are essentially negated at very small timestep sizes.  While this finding likely does not effect results using current production timestep sizes, it needs to be considered as spatial resolution is increased (effectively reducing production timestep size).  Furthermore, the energy conservation error from the vertical remap approach is positive, indicating that the remap process adds energy to the system.  The authors therefore recommend investigating possible improvements to these approaches in HOMME-NH.  One way to address the hyperviscosity error is to include the application of hyperviscocity in the explicit portion of the ARK IMEX scheme instead of applying it in a split fashion.  Meanwhile, the vertical remap approach might be swapped out for one that discretely conserves energy or, if that is not feasible, for one that dissipates a small amount of energy.  This might allow for a smaller hyperviscosity coefficient $\nu$, because the results in Figure \ref{FIG:baroclinic_energy_error} suggest the coefficient is currently tuned to counteract the energy added by the vertical remap.

\section{Conclusions} \label{SEC:conclusions}
  The non-hydrostatic atmosphere model, presented in Section \ref{SEC:nonhydrostatic}, provides the interesting challenge of handling acoustic wave propagation.  This work seeks to find the time integration scheme to be implemented in the HOMME-NH dynamical core that most efficiently addresses this challenge.  A class of time integration schemes known as implicit-explicit additive Runge-Kutta (ARK IMEX) methods are selected for evaluation due to their ability to treat troublesome model terms implicitly in time while treating others in an  explicit manner, all without losing access to higher-order accuracy.  27 different ARK IMEX methods, listed in Section \ref{SEC:numerical}, are evaluated on efficiency in producing an accurate solution, ability to take the corresponding hydrostatic timestep, and sensitivity to horizontal grid resolution.  Instead of a single ideal method, the authors recommend 5 ARK IMEX methods be implemented in HOMME-NH and suggest that improvements to the hyperviscosity and vertical remap approaches will significantly improve time integration performance.

  In Section \ref{SEC:evaluation}, the DCMIP 2012 Non-orographic Gravity Wave and Baroclinic Instability tests are used to evaluate the ARK IMEX methods.  A convergence test first verifies that the ARK IMEX methods, implemented with ARKode, attain the expected convergence orders in the absence of hyperviscosity.  When hyperviscosity is applied in split fashion, however, the test reveals that the error from the hyperviscosity approach in HOMME-NH can dominate the overall error for small enough timestep sizes.  This finding was verified with global energy conservation in the Baroclinic Instability test, where the errors from the hyperviscosity and vertical remap approaches were compared to the error from the ARK IMEX method itself.  It was found that for small timestep sizes, both the hyperviscosity and vertical remap errors are orders of magnitude larger than the ARK IMEX method, with the vertical remap additionally causing an increase in the global energy of the system.  These small timestep sizes are currently below the production timestep size, but increases in spatial resolution will continue to decrease that production timestep size.  For this reason, the authors recommend possibly including the hyperviscosity application in the explicit portion of the ARK IMEX method, to reduce the hyperviscosity error and modifying the vertical remap approach to either discretely conserve energy or be slightly dissipative.

  For performance at the current production timestep size, the 27 ARK IMEX methods are evaluated on how efficiently an accurate 15-day solution can be produced in the Baroclinic Instability test, which requires both hyperviscosity and vertical remap.  An accuracy criterion is defined using the surface pressure, where the error tolerance is defined as the difference between the solution using a production timestep size and a reference solution using a timestep size that is $30$ times smaller.  In addition to producing a solution of acceptable accuracy in an efficient manner, the ARK IMEX methods are evaluated on the ability to produce a solution, regardless of accuracy, using the hydrostatic timestep of $\Delta t=300$s.  Finally, the sensitivity of the methods to shrinking horizontal grid-cell lengths was investigated using planets of reduced radii.  The methods that excelled in each metric are summarized in Section \ref{SEC:recommendations}, with the authors recommending that the ARK324, DBM453, IMKG242b, IMKG243a, and IMKG252b methods be implemented in HOMME-NH.  The authors encourage future work using the metrics described herein to test newly developed integration techniques, ARK IMEX or otherwise, on more complex test problems (full planet nonhydrostatic dynamics, orography, etc) using both HOMME-NH and other dynamical cores.

\appendix

\section{Custom IMEX ARK method} \label{APP:dbm}
The DBM453 method is defined by explicit and implicit Butcher tables with shared abcissae $c_E = c_I =
( 0, 0.1030620881159184, 0.72139131281753662, 1.28181117351981733, 1)$ and shared root nodes $b_E = b_I = ( 0.87795339639076672, -0.72692641526151549,$ \\ \noindent $0.7520413715737272, -0.22898029400415090, 0.32591194130117247)$.  The explicit coefficients are
\small
\[
   A_E = \begin{bmatrix} 0& 0& 0& 0& 0\\
         0.10306208811591838& 0& 0& 0& 0\\
        -0.94124866143519894& 1.6626399742527356& 0& 0& 0\\
        -1.3670975201437765& 1.3815852911016873& 1.2673234025619065& 0& 0\\
        -0.81287582068772448& 0.81223739060505738& 0.90644429603699305& 0.094194134045674111& 0\end{bmatrix}
\]
\normalsize
and the implicit coefficients are
\small
\[
  A_I = \begin{bmatrix} 0& 0& 0& 0& 0\\
        -0.22284985318525410& \gamma& 0& 0& 0\\
        -0.46801347074080545& 0.86349284225716961& \gamma& 0& 0\\
        -0.46509906651927421& 0.81063103116959553& 0.61036726756832357& \gamma& 0\\
        0.87795339639076675& -0.72692641526151547& 0.75204137157372720& -0.22898029400415088& \gamma \end{bmatrix}
\]
\normalsize
where $\gamma = 0.32591194130117247$.  The method was designed to satisfy the following properties:
\begin{itemize}
\item[(a)] 5 stages, with 4 implicit solves per step, with the same coefficient $\gamma$ used in each implicit solve,
\item[(b)] third order accuracy for both explicit and implicit methods, with third order coupling,
\item[(c)] the explicit method has provably maximal stability along the imaginary axis for a 5-stage method (see \cite{kinnmark_one_1984}), and
\item[(d)] the implicit method is both $A$- and $L$-stable.
\end{itemize}
Of these, the authors felt that property (c) was of fundamental importance for this application, due to the fact that the eigenvalues of the explicit portion of the model, $\mathbf{f}^E(\mathbf{q})$ lie on the imaginary axis, and frequently serve to limit the maximum stable step size for IMEX methods.

\section{IMEX ARK method properties} \label{APP:ARK_properties}
Table \ref{tab:ark_properties} provides a variety of theoretical properties of each of the ARK methods used in this paper.  The categories match those from \citeA<>[Table A1]{gardner_implicitexplicit_2018}:
\begin{itemize}
    \item number of implicit solves per step ($\mathbf{f}^I$ column);
    \item number of explicit stages per step ($\mathbf{f}^E$ column);
    \item theoretical order of accuracy for the explicit, diagonally-implicit, and overall additive Runge--Kutta methods (Order E, I and A);
    \item theoretical order of accuracy for internal stages (relevant for order reduction on stiff problems) for the ERK, DIRK and ARK stages (Stage order E, I and A);
    \item whether the DIRK method is A-, L- and/or B-stable (Stabilty A, L and B);
    \item whether the DIRK and ERK methods are stiffly accurate, i.e., the last row of $A^I$ equals $b^I$, and last row of $A^E$ equals $b^E$ (S.A. DIRK and ERK);
    \item whether the DIRK and ERK methods share the same solution weights (i.e., $b^I=b^E$) and thus the method will preserve linear invariants to machine precision ($b$);
    \item whether the DIRK and ERK methods share the same abcissae (i.e. $c^I=c^E$) and thus $\mathbf{f}^I$ and $\mathbf{f}^E$ are evaluated at the same stage times ($c$);
    \item the maximum stable explicit step along the imaginary axis i.e., the largest $y_{max}\ge 0$ such that the ERK method is stable for all $\lambda=iy, 0\le y\le y_{max}$ (Max exp).
\end{itemize}
These are assessed using the same mechanisms as outlined in \citeA{gardner_implicitexplicit_2018}.

% Table of ARK properties
% \ding{51} is a check-mark
% \ding{55} is an x-mark
\begin{table*}[tbh]
\caption{Properties for each of the ARK methods used in this paper. The column headings are described in the above text. \label{tab:ark_properties}}
\begin{tabular}{lccccc|ccc|cccccccc}
\toprule
 \multirow{2}{*}{Method}&
 \multirow{2}{*}{$\mathbf{f}^I$}&
 \multirow{2}{*}{$\textbf{f}^E$}&
 \multicolumn{3}{c}{Order}&
 \multicolumn{3}{c}{Stage order}&
 \multicolumn{3}{c}{Stability}&
 S.A.&
 S.A.&
 \multirow{2}{*}{$b$}&
 \multirow{2}{*}{$c$}&
 \multirow{2}{*}{Max Exp}\\
 & & & E & I & A & E & I & A & A & L & B & DIRK & ERK & & & \\
 \midrule
 ARS222  & 2 & 3 & 2 & 2 & 2 & 1 & 1 & 1 & \ding{51} & \ding{51} & \ding{55} & \ding{51} & \ding{51} & \ding{55} & \ding{51} &  $\sim\!\! 0.00$ \\
 ARS232  & 2 & 3 & 2 & 2 & 2 & 1 & 1 & 1 & \ding{51} & \ding{51} & \ding{55} & \ding{51} & \ding{55} & \ding{51} & \ding{51} &  $\sim\!\! 1.73$ \\
  GSA222  & 2 & 3 & 2 & 2 & 2 & 1 & 1 & 1 & \ding{51} & \ding{51} & \ding{55} & \ding{51} & \ding{51} & \ding{55} & \ding{51} &  $\sim\!\! 0.00$ \\
 SSP2232 & 2 & 3 & 2 & 2 & 2 & 1 & 2 & 0 & \ding{51} & \ding{55} & \ding{55} & \ding{51} & \ding{55} & \ding{51} & \ding{55} &  $\sim\!\! 1.73$ \\
 IMKG232(a,b)& 2 & 3 & 2 & 2 & 2 & 1 & 1 & 0 & \ding{51} & \ding{51} & \ding{55} & \ding{51} & \ding{51} & \ding{51} & \ding{55} &  $\sim\!\! 2.00$ \\
 IMKG242(a,b)& 2 & 4 & 2 & 2 & 2 & 1 & 1 & 0 & \ding{51} & \ding{51} & \ding{55} & \ding{51} & \ding{51} & \ding{51} & \ding{55} &  $\sim\!\! 2.83$ \\
 IMKG252(a,b)& 2 & 5 & 2 & 2 & 2 & 1 & 1 & 0 & \ding{51} & \ding{51} & \ding{55} & \ding{51} & \ding{51} & \ding{51} & \ding{55} &  $\sim\!\! 4.00$ \\
 IMKG243a& 3 & 4 & 2 & 2 & 2 & 1 & 1 & 0 & \ding{51} & \ding{51} & \ding{55} & \ding{51} & \ding{51} & \ding{51} & \ding{55} &  $\sim\!\! 2.83$ \\
 IMKG253(a,b)& 3 & 5 & 2 & 2 & 2 & 1 & 1 & 0 & \ding{51} & \ding{51} & \ding{55} & \ding{51} & \ding{51} & \ding{51} & \ding{55} &  $\sim\!\! 4.00$ \\
 IMKG254(a,b)& 4 & 5 & 2 & 2 & 2 & 1 & 1 & 0 & \ding{55} & \ding{55} & \ding{55} & \ding{51} & \ding{51} & \ding{51} & \ding{55} &  $\sim\!\! 4.00$ \\
 IMKG254c& 4 & 5 & 2 & 2 & 2 & 1 & 1 & 0 & \ding{51} & \ding{51} & \ding{55} & \ding{51} & \ding{51} & \ding{51} & \ding{55} &  $\sim\!\! 4.00$ \\
 \midrule
 ARS233  & 2 & 3 & 3 & 3 & 3 & 1 & 1 & 1 & \ding{51} & \ding{55} & \ding{51} & \ding{55} & \ding{55} & \ding{51} & \ding{51} &  $\sim\!\! 1.73$ \\
 SSP3333b& 2 & 3 & 3 & 3 & 3 & 1 & 1 & 1 & \ding{51} & \ding{55} & \ding{55} & \ding{55} & \ding{55} & \ding{51} & \ding{51} &  $\sim\!\! 1.73$ \\
 SSP3333c& 2 & 3 & 3 & 3 & 3 & 1 & 1 & 1 & \ding{51} & \ding{55} & \ding{55} & \ding{55} & \ding{55} & \ding{51} & \ding{51} &  $\sim\!\! 1.73$ \\
 ARK324  & 3 & 4 & 3 & 3 & 3 & 1 & 2 & 1 & \ding{51} & \ding{51} & \ding{55} & \ding{51} & \ding{55} & \ding{51} & \ding{51} &  $\sim\!\! 2.48$ \\
 ARS343  & 3 & 4 & 3 & 3 & 3 & 1 & 1 & 1 & \ding{51} & \ding{51} & \ding{55} & \ding{51} & \ding{55} & \ding{51} & \ding{51} &  $\sim\!\! 2.83$ \\
 IMKG343a& 3 & 4 & 3 & 3 & 3 & 1 & 1 & 0 & \ding{55} & \ding{55} & \ding{55} & \ding{51} & \ding{51} & \ding{51} & \ding{55} &  $\sim\!\! 2.83$ \\
 ARS443  & 4 & 4 & 3 & 3 & 3 & 1 & 1 & 1 & \ding{51} & \ding{51} & \ding{55} & \ding{51} & \ding{51} & \ding{55} & \ding{51} &  $\sim\!\! 1.57$ \\
 DBM453  & 4 & 5 & 3 & 3 & 3 & 1 & 1 & 1 & \ding{51} & \ding{51} & \ding{55} & \ding{51} & \ding{55} & \ding{51} & \ding{51} &  $\sim\!\! 3.87$ \\
 ARK436  & 5 & 6 & 4 & 4 & 4 & 1 & 2 & 1 & \ding{51} & \ding{51} & \ding{55} & \ding{51} & \ding{55} & \ding{51} & \ding{51} &  $\sim\!\! 4.00$ \\
 ARK437  & 6 & 7 & 4 & 4 & 4 & 1 & 2 & 1 & \ding{51} & \ding{51} & \ding{55} & \ding{51} & \ding{55} & \ding{51} & \ding{51} &  $\sim\!\! 4.70$ \\
 ARK548  & 7 & 8 & 5 & 5 & 5 & 1 & 2 & 1 & \ding{51} & \ding{51} & \ding{55} & \ding{51} & \ding{55} & \ding{51} & \ding{51} &  $\sim\!\! 0.02$ \\
\bottomrule
\end{tabular}
%\belowtable{} % Table Footnotes
\end{table*}

\acknowledgments
The authors first wish to acknowledge and express gratitude to David J. Gardner for his help in interfacing the ARKode package with HOMME-NH.

The HOMME-NH and E3SM code, namelist files, post-processing scripts, post-processed data, and plotting scripts are all available in a Zenodo archive \cite{vogl_code_2019}.  That archive is a snapshot of \texttt{homme-nh\_ARK\_IMEX} branch of the GitHub repository at \url{https://github.com/cjvogl/E3SM}, which is a fork of the E3SM GitHub repository \cite{e3sm_project_energy_2018} with the addition of content for reproducing the results presented herein.  Further description is in \texttt{README.txt} located in the \texttt{components/homme/scripts\_for\_paper} directory of the archive.  A second Zenodo archive \cite{vogl_stability_2019} provides an archive of plots showing the stability regions for each ARK IMEX method.

This work was supported in part by the U.S. Department of Energy, Office of Science, Office of Advanced Scientific Computing Research and Office of Biological and Environmental Research through the partnership project, ``A Non-hydrostatic Variable Resolution Atmospheric Model in ACME.''

This work was performed under the auspices of the U.S. Department of Energy by Lawrence Livermore National Laboratory under contract DE-AC52-07NA27344. Lawrence Livermore National Security, LLC.  LLNL-JRNL-770788-DRAFT.

Sandia National Laboratories is a multi-mission laboratory managed and operated by National Technology and Engineering Solutions of Sandia, LLC., a wholly owned subsidiary of Honeywell International, Inc., for the U.S. Department of Energy’s National Nuclear Security Administration under contract DE-NA0003525.  This paper describes objective technical results and analysis. Any subjective views or opinions that might be expressed in the paper do not necessarily represent the views of the U.S. Department of Energy or the United States Government.

This document was prepared as an account of work sponsored by an agency of the United States government. Neither the United States government nor Lawrence Livermore National Security, LLC, nor any of their employees makes any warranty, expressed or implied, or assumes any legal liability or responsibility for the accuracy, completeness, or usefulness of any information, apparatus, product, or process disclosed, or represents that its use would not infringe privately owned rights. Reference herein to any specific commercial product, process, or service by trade name, trademark, manufacturer, or otherwise does not necessarily constitute or imply its endorsement, recommendation, or favoring by the United States government or Lawrence Livermore National Security, LLC. The views and opinions of authors expressed herein do not necessarily state or reflect those of the United States government or Lawrence Livermore National Security, LLC, and shall not be used for advertising or product endorsement purposes.

%% ------------------------------------------------------------------------ %%
%% References and Citations

%%%%%%%%%%%%%%%%%%%%%%%%%%%%%%%%%%%%%%%%%%%%%%%
% BibTeX is preferred:
%
% \bibliography{<name of your .bib file>}

\begin{thebibliography}{}

\bibitem [\protect \citeauthoryear {%
Arakawa%
\ \BBA {} Konor%
}{%
Arakawa%
\ \BBA {} Konor%
}{%
{\protect \APACyear {2009}}%
}]{%
arakawa_unification_2009}
\APACinsertmetastar {%
arakawa_unification_2009}%
\begin{APACrefauthors}%
Arakawa, A.%
\BCBT {}\ \BBA {} Konor, C\BPBI S.%
\end{APACrefauthors}%
\unskip\
\newblock
\APACrefYearMonthDay{2009}{}{}.
\newblock
{\BBOQ}\APACrefatitle {Unification of the anelastic and quasi-hydrostatic
  systems of equations} {Unification of the anelastic and quasi-hydrostatic
  systems of equations}.{\BBCQ}
\newblock
\APACjournalVolNumPages{Mon. Wea. Rev.}{137}{2}{710--726}.
\newblock
\begin{APACrefDOI} \doi{10.1175/2008MWR2520.1} \end{APACrefDOI}
\PrintBackRefs{\CurrentBib}

\bibitem [\protect \citeauthoryear {%
Araujo%
, Murua%
\BCBL {}\ \BBA {} Sanz-Serna%
}{%
Araujo%
\ \protect \BOthers {.}}{%
{\protect \APACyear {1997}}%
}]{%
araujo_symplectic_1997}
\APACinsertmetastar {%
araujo_symplectic_1997}%
\begin{APACrefauthors}%
Araujo, A.%
, Murua, A.%
\BCBL {}\ \BBA {} Sanz-Serna, J.%
\end{APACrefauthors}%
\unskip\
\newblock
\APACrefYearMonthDay{1997}{}{}.
\newblock
{\BBOQ}\APACrefatitle {Symplectic methods based on decompositions} {Symplectic
  methods based on decompositions}.{\BBCQ}
\newblock
\APACjournalVolNumPages{SIAM J. Numer. Anal.}{34}{5}{1926--1947}.
\newblock
\begin{APACrefDOI} \doi{10.1137/S0036142995292128} \end{APACrefDOI}
\PrintBackRefs{\CurrentBib}

\bibitem [\protect \citeauthoryear {%
Ascher%
, Ruuth%
\BCBL {}\ \BBA {} Spiteri%
}{%
Ascher%
\ \protect \BOthers {.}}{%
{\protect \APACyear {1997}}%
}]{%
ascher_implicit-explicit_1997}
\APACinsertmetastar {%
ascher_implicit-explicit_1997}%
\begin{APACrefauthors}%
Ascher, U\BPBI M.%
, Ruuth, S\BPBI J.%
\BCBL {}\ \BBA {} Spiteri, R\BPBI J.%
\end{APACrefauthors}%
\unskip\
\newblock
\APACrefYearMonthDay{1997}{}{}.
\newblock
{\BBOQ}\APACrefatitle {Implicit-explicit {Runge--Kutta} methods for
  time-dependent partial differential equations} {Implicit-explicit
  {Runge--Kutta} methods for time-dependent partial differential
  equations}.{\BBCQ}
\newblock
\APACjournalVolNumPages{Appl. Numer. Math.}{25}{2-3}{151--167}.
\newblock
\begin{APACrefDOI} \doi{10.1016/S0168-9274(97)00056-1} \end{APACrefDOI}
\PrintBackRefs{\CurrentBib}

\bibitem [\protect \citeauthoryear {%
Boscarino%
, Russo%
\BCBL {}\ \BBA {} Scandurra%
}{%
Boscarino%
\ \protect \BOthers {.}}{%
{\protect \APACyear {2018}}%
}]{%
boscarino_all_2018}
\APACinsertmetastar {%
boscarino_all_2018}%
\begin{APACrefauthors}%
Boscarino, S.%
, Russo, G.%
\BCBL {}\ \BBA {} Scandurra, L.%
\end{APACrefauthors}%
\unskip\
\newblock
\APACrefYearMonthDay{2018}{{\APACmonth{11}}}{}.
\newblock
{\BBOQ}\APACrefatitle {All {Mach} {Number} {Second} {Order} {Semi}-implicit
  {Scheme} for the {Euler} {Equations} of {Gas} {Dynamics}} {All {Mach}
  {Number} {Second} {Order} {Semi}-implicit {Scheme} for the {Euler}
  {Equations} of {Gas} {Dynamics}}.{\BBCQ}
\newblock
\APACjournalVolNumPages{J. Sci. Comput.}{77}{2}{850--884}.
\newblock
\begin{APACrefDOI} \doi{10.1007/s10915-018-0731-9} \end{APACrefDOI}
\PrintBackRefs{\CurrentBib}

\bibitem [\protect \citeauthoryear {%
Conde%
, Gottlieb%
, Grant%
\BCBL {}\ \BBA {} Shadid%
}{%
Conde%
\ \protect \BOthers {.}}{%
{\protect \APACyear {2017}}%
}]{%
conde_implicit_2017}
\APACinsertmetastar {%
conde_implicit_2017}%
\begin{APACrefauthors}%
Conde, S.%
, Gottlieb, S.%
, Grant, Z\BPBI J.%
\BCBL {}\ \BBA {} Shadid, J\BPBI N.%
\end{APACrefauthors}%
\unskip\
\newblock
\APACrefYearMonthDay{2017}{}{}.
\newblock
{\BBOQ}\APACrefatitle {Implicit and implicit-explicit strong stability
  preserving {Runge-–Kutta} methods with high linear order} {Implicit and
  implicit-explicit strong stability preserving {Runge-–Kutta} methods with
  high linear order}.{\BBCQ}
\newblock
\APACjournalVolNumPages{J. Sci. Comput.}{73}{2-3}{667--690}.
\newblock
\begin{APACrefDOI} \doi{10.1007/s10915-017-0560-2} \end{APACrefDOI}
\PrintBackRefs{\CurrentBib}

\bibitem [\protect \citeauthoryear {%
Davies%
, Staniforth%
, Wood%
\BCBL {}\ \BBA {} Thuburn%
}{%
Davies%
\ \protect \BOthers {.}}{%
{\protect \APACyear {2003}}%
}]{%
davies_validity_2003}
\APACinsertmetastar {%
davies_validity_2003}%
\begin{APACrefauthors}%
Davies, T.%
, Staniforth, A.%
, Wood, N.%
\BCBL {}\ \BBA {} Thuburn, J.%
\end{APACrefauthors}%
\unskip\
\newblock
\APACrefYearMonthDay{2003}{}{}.
\newblock
{\BBOQ}\APACrefatitle {Validity of anelastic and other equation sets as
  inferred from normal-mode analysis} {Validity of anelastic and other equation
  sets as inferred from normal-mode analysis}.{\BBCQ}
\newblock
\APACjournalVolNumPages{Quart. J. Roy. Meteorol. Soc.}{129}{593}{2761--2775}.
\newblock
\begin{APACrefDOI} \doi{10.1256/qj.02.1951} \end{APACrefDOI}
\PrintBackRefs{\CurrentBib}

\bibitem [\protect \citeauthoryear {%
Dennis%
\ \protect \BOthers {.}}{%
Dennis%
\ \protect \BOthers {.}}{%
{\protect \APACyear {2012}}%
}]{%
dennis_cam-se:_2012}
\APACinsertmetastar {%
dennis_cam-se:_2012}%
\begin{APACrefauthors}%
Dennis, J\BPBI M.%
, Edwards, J.%
, Evans, K\BPBI J.%
, Guba, O.%
, Lauritzen, P\BPBI H.%
, Mirin, A\BPBI A.%
\BDBL {}Worley, P\BPBI H.%
\end{APACrefauthors}%
\unskip\
\newblock
\APACrefYearMonthDay{2012}{}{}.
\newblock
{\BBOQ}\APACrefatitle {{CAM}-{SE}: {A} scalable spectral element dynamical core
  for the {Community} {Atmosphere} {Model}} {{CAM}-{SE}: {A} scalable spectral
  element dynamical core for the {Community} {Atmosphere} {Model}}.{\BBCQ}
\newblock
\APACjournalVolNumPages{Internat. J. High Perf. Comput. Appl.}{26}{1}{74--89}.
\newblock
\begin{APACrefDOI} \doi{10.1177/1094342011428142} \end{APACrefDOI}
\PrintBackRefs{\CurrentBib}

\bibitem [\protect \citeauthoryear {%
Durran%
}{%
Durran%
}{%
{\protect \APACyear {1989}}%
}]{%
durran_improving_1989}
\APACinsertmetastar {%
durran_improving_1989}%
\begin{APACrefauthors}%
Durran, D\BPBI R.%
\end{APACrefauthors}%
\unskip\
\newblock
\APACrefYearMonthDay{1989}{}{}.
\newblock
{\BBOQ}\APACrefatitle {Improving the anelastic approximation} {Improving the
  anelastic approximation}.{\BBCQ}
\newblock
\APACjournalVolNumPages{J. Atmos. Sci.}{46}{11}{1453--1461}.
\newblock
\begin{APACrefDOI} \doi{10.1175/1520-0469(1989)046<1453:ITAA>2.0.CO;2}
  \end{APACrefDOI}
\PrintBackRefs{\CurrentBib}

\bibitem [\protect \citeauthoryear {%
{E3SM Project}%
}{%
{E3SM Project}%
}{%
{\protect \APACyear {2018}}%
}]{%
e3sm_project_energy_2018}
\APACinsertmetastar {%
e3sm_project_energy_2018}%
\begin{APACrefauthors}%
{E3SM Project}.%
\end{APACrefauthors}%
\unskip\
\newblock
\APACrefYearMonthDay{2018}{}{}.
\newblock
\APACrefbtitle {Energy {Exascale} {Earth} {System} {Model} ({E}3{SM}).} {Energy
  {Exascale} {Earth} {System} {Model} ({E}3{SM}).}
\newblock
\APACrefnote{\url{https://dx.doi.org/10.11578/E3SM/dc.20180418.36}}
\newblock
\begin{APACrefDOI} \doi{10.11578/E3SM/dc.20180418.36} \end{APACrefDOI}
\PrintBackRefs{\CurrentBib}

\bibitem [\protect \citeauthoryear {%
Evans%
\ \protect \BOthers {.}}{%
Evans%
\ \protect \BOthers {.}}{%
{\protect \APACyear {2017}}%
}]{%
evans_performance_2017}
\APACinsertmetastar {%
evans_performance_2017}%
\begin{APACrefauthors}%
Evans, K\BPBI J.%
, Archibald, R\BPBI K.%
, Gardner, D\BPBI J.%
, Norman, M\BPBI R.%
, Taylor, M\BPBI A.%
, Woodward, C\BPBI S.%
\BCBL {}\ \BBA {} Worley, P\BPBI H.%
\end{APACrefauthors}%
\unskip\
\newblock
\APACrefYearMonthDay{2017}{}{}.
\newblock
{\BBOQ}\APACrefatitle {Performance analysis of fully explicit and fully
  implicit solvers within a spectral element shallow-water atmosphere model}
  {Performance analysis of fully explicit and fully implicit solvers within a
  spectral element shallow-water atmosphere model}.{\BBCQ}
\newblock
\APACjournalVolNumPages{Internat. J. High Perf. Comput.
  Appl.}{}{}{1094342017736373}.
\newblock
\begin{APACrefDOI} \doi{10.1177/1094342017736373} \end{APACrefDOI}
\PrintBackRefs{\CurrentBib}

\bibitem [\protect \citeauthoryear {%
Gardner%
\ \protect \BOthers {.}}{%
Gardner%
\ \protect \BOthers {.}}{%
{\protect \APACyear {2018}}%
}]{%
gardner_implicitexplicit_2018}
\APACinsertmetastar {%
gardner_implicitexplicit_2018}%
\begin{APACrefauthors}%
Gardner, D\BPBI J.%
, Guerra, J\BPBI E.%
, Hamon, F\BPBI P.%
, Reynolds, D\BPBI R.%
, Ullrich, P\BPBI A.%
\BCBL {}\ \BBA {} Woodward, C\BPBI S.%
\end{APACrefauthors}%
\unskip\
\newblock
\APACrefYearMonthDay{2018}{}{}.
\newblock
{\BBOQ}\APACrefatitle {Implicit-explicit ({IMEX}) {Runge-–Kutta} methods for
  non-hydrostatic atmospheric models} {Implicit-explicit ({IMEX})
  {Runge-–Kutta} methods for non-hydrostatic atmospheric models}.{\BBCQ}
\newblock
\APACjournalVolNumPages{Geophys. Model Dev.}{11}{4}{1497--1515}.
\newblock
\begin{APACrefDOI} \doi{https://doi.org/10.5194/gmd-11-1497-2018}
  \end{APACrefDOI}
\PrintBackRefs{\CurrentBib}

\bibitem [\protect \citeauthoryear {%
Giraldo%
, Kelly%
\BCBL {}\ \BBA {} Constantinescu%
}{%
Giraldo%
\ \protect \BOthers {.}}{%
{\protect \APACyear {2013}}%
}]{%
giraldo_implicit-explicit_2013}
\APACinsertmetastar {%
giraldo_implicit-explicit_2013}%
\begin{APACrefauthors}%
Giraldo, F.%
, Kelly, J.%
\BCBL {}\ \BBA {} Constantinescu, E.%
\end{APACrefauthors}%
\unskip\
\newblock
\APACrefYearMonthDay{2013}{}{}.
\newblock
{\BBOQ}\APACrefatitle {Implicit-explicit formulations of a three-dimensional
  nonhydrostatic unified model of the atmosphere ({NUMA})} {Implicit-explicit
  formulations of a three-dimensional nonhydrostatic unified model of the
  atmosphere ({NUMA})}.{\BBCQ}
\newblock
\APACjournalVolNumPages{SIAM J. Sci. Comput.}{35}{5}{B1162--B1194}.
\newblock
\begin{APACrefDOI} \doi{10.1137/120876034} \end{APACrefDOI}
\PrintBackRefs{\CurrentBib}

\bibitem [\protect \citeauthoryear {%
Guerra%
\ \BBA {} Ullrich%
}{%
Guerra%
\ \BBA {} Ullrich%
}{%
{\protect \APACyear {2016}}%
}]{%
guerra_high-order_2016}
\APACinsertmetastar {%
guerra_high-order_2016}%
\begin{APACrefauthors}%
Guerra, J\BPBI E.%
\BCBT {}\ \BBA {} Ullrich, P\BPBI A.%
\end{APACrefauthors}%
\unskip\
\newblock
\APACrefYearMonthDay{2016}{}{}.
\newblock
{\BBOQ}\APACrefatitle {A high-order staggered finite-element vertical
  discretization for non-hydrostatic atmospheric models} {A high-order
  staggered finite-element vertical discretization for non-hydrostatic
  atmospheric models}.{\BBCQ}
\newblock
\APACjournalVolNumPages{Geophys. Model Dev.}{9}{5}{2007--2029}.
\newblock
\begin{APACrefDOI} \doi{https://doi.org/10.5194/gmd-9-2007-2016}
  \end{APACrefDOI}
\PrintBackRefs{\CurrentBib}

\bibitem [\protect \citeauthoryear {%
Harris%
, Lin%
\BCBL {}\ \BBA {} Tu%
}{%
Harris%
\ \protect \BOthers {.}}{%
{\protect \APACyear {2016}}%
}]{%
harris_high-resolution_2016}
\APACinsertmetastar {%
harris_high-resolution_2016}%
\begin{APACrefauthors}%
Harris, L\BPBI M.%
, Lin, S\BHBI J.%
\BCBL {}\ \BBA {} Tu, C.%
\end{APACrefauthors}%
\unskip\
\newblock
\APACrefYearMonthDay{2016}{}{}.
\newblock
{\BBOQ}\APACrefatitle {High-resolution climate simulations using {GFDL} {HiRAM}
  with a stretched global grid} {High-resolution climate simulations using
  {GFDL} {HiRAM} with a stretched global grid}.{\BBCQ}
\newblock
\APACjournalVolNumPages{J. Climate}{29}{11}{4293--4314}.
\newblock
\begin{APACrefDOI} \doi{10.1175/JCLI-D-15-0389.1} \end{APACrefDOI}
\PrintBackRefs{\CurrentBib}

\bibitem [\protect \citeauthoryear {%
Hindmarsh%
\ \protect \BOthers {.}}{%
Hindmarsh%
\ \protect \BOthers {.}}{%
{\protect \APACyear {2005}}%
}]{%
hindmarsh_sundials:_2005}
\APACinsertmetastar {%
hindmarsh_sundials:_2005}%
\begin{APACrefauthors}%
Hindmarsh, B\BPBI P\BPBI N., Alan~C.%
, Grant, K\BPBI E.%
, Lee, S\BPBI L.%
, Serban, R.%
, Shumaker, D\BPBI E.%
\BCBL {}\ \BBA {} Woodward, C\BPBI S.%
\end{APACrefauthors}%
\unskip\
\newblock
\APACrefYearMonthDay{2005}{}{}.
\newblock
{\BBOQ}\APACrefatitle {{SUNDIALS}: {Suite} of {Nonlinear} and
  {Differential}/{Algebraic} {Equation} {Solvers}} {{SUNDIALS}: {Suite} of
  {Nonlinear} and {Differential}/{Algebraic} {Equation} {Solvers}}.{\BBCQ}
\newblock
\APACjournalVolNumPages{ACM Trans. Math. Software}{31}{3}{363--396}.
\newblock
\begin{APACrefDOI} \doi{10.1145/1089014.1089020} \end{APACrefDOI}
\PrintBackRefs{\CurrentBib}

\bibitem [\protect \citeauthoryear {%
Huang%
, Rhoades%
, Ullrich%
\BCBL {}\ \BBA {} Zarzycki%
}{%
Huang%
\ \protect \BOthers {.}}{%
{\protect \APACyear {2016}}%
}]{%
huang_evaluation_2016}
\APACinsertmetastar {%
huang_evaluation_2016}%
\begin{APACrefauthors}%
Huang, X.%
, Rhoades, A\BPBI M.%
, Ullrich, P\BPBI A.%
\BCBL {}\ \BBA {} Zarzycki, C\BPBI M.%
\end{APACrefauthors}%
\unskip\
\newblock
\APACrefYearMonthDay{2016}{}{}.
\newblock
{\BBOQ}\APACrefatitle {An evaluation of the variable-resolution {CESM} for
  modeling {California}'s climate} {An evaluation of the variable-resolution
  {CESM} for modeling {California}'s climate}.{\BBCQ}
\newblock
\APACjournalVolNumPages{J. Adv. Mod. Earth Sys.}{8}{1}{345--369}.
\newblock
\begin{APACrefDOI} \doi{10.1002/2015MS000559} \end{APACrefDOI}
\PrintBackRefs{\CurrentBib}

\bibitem [\protect \citeauthoryear {%
Jablonowski%
\ \BBA {} Williamson%
}{%
Jablonowski%
\ \BBA {} Williamson%
}{%
{\protect \APACyear {2006}}%
}]{%
jablonowski_baroclinic_2006}
\APACinsertmetastar {%
jablonowski_baroclinic_2006}%
\begin{APACrefauthors}%
Jablonowski, C.%
\BCBT {}\ \BBA {} Williamson, D\BPBI L.%
\end{APACrefauthors}%
\unskip\
\newblock
\APACrefYearMonthDay{2006}{}{}.
\newblock
{\BBOQ}\APACrefatitle {A baroclinic instability test case for atmospheric model
  dynamical cores} {A baroclinic instability test case for atmospheric model
  dynamical cores}.{\BBCQ}
\newblock
\APACjournalVolNumPages{Quart. J. Roy. Meteorol. Soc.}{132}{621C}{2943--2975}.
\newblock
\begin{APACrefDOI} \doi{10.1256/qj.06.12} \end{APACrefDOI}
\PrintBackRefs{\CurrentBib}

\bibitem [\protect \citeauthoryear {%
Kennedy%
\ \BBA {} Carpenter%
}{%
Kennedy%
\ \BBA {} Carpenter%
}{%
{\protect \APACyear {2003}}%
}]{%
kennedy_additive_2001}
\APACinsertmetastar {%
kennedy_additive_2001}%
\begin{APACrefauthors}%
Kennedy, C\BPBI A.%
\BCBT {}\ \BBA {} Carpenter, M\BPBI H.%
\end{APACrefauthors}%
\unskip\
\newblock
\APACrefYearMonthDay{2003}{}{}.
\newblock
{\BBOQ}\APACrefatitle {Additive {Runge-–Kutta} schemes for
  convection–diffusion–reaction equations} {Additive {Runge-–Kutta}
  schemes for convection–diffusion–reaction equations}.{\BBCQ}
\newblock
\APACjournalVolNumPages{Applied Numerical Mathematics}{44}{1}{139 - 181}.
\newblock
\begin{APACrefDOI} \doi{10.1016/S0168-9274(02)00138-1} \end{APACrefDOI}
\PrintBackRefs{\CurrentBib}

\bibitem [\protect \citeauthoryear {%
Kennedy%
\ \BBA {} Carpenter%
}{%
Kennedy%
\ \BBA {} Carpenter%
}{%
{\protect \APACyear {2019}}%
}]{%
kennedy_higher-order_2019}
\APACinsertmetastar {%
kennedy_higher-order_2019}%
\begin{APACrefauthors}%
Kennedy, C\BPBI A.%
\BCBT {}\ \BBA {} Carpenter, M\BPBI H.%
\end{APACrefauthors}%
\unskip\
\newblock
\APACrefYearMonthDay{2019}{{\APACmonth{02}}}{}.
\newblock
{\BBOQ}\APACrefatitle {Higher-order additive {Runge}–{Kutta} schemes for
  ordinary differential equations} {Higher-order additive {Runge}–{Kutta}
  schemes for ordinary differential equations}.{\BBCQ}
\newblock
\APACjournalVolNumPages{Appl. Numer. Math.}{136}{}{183--205}.
\newblock
\begin{APACrefDOI} \doi{10.1016/j.apnum.2018.10.007} \end{APACrefDOI}
\PrintBackRefs{\CurrentBib}

\bibitem [\protect \citeauthoryear {%
Kinnmark%
\ \BBA {} Gray%
}{%
Kinnmark%
\ \BBA {} Gray%
}{%
{\protect \APACyear {1984}}%
}]{%
kinnmark_one_1984}
\APACinsertmetastar {%
kinnmark_one_1984}%
\begin{APACrefauthors}%
Kinnmark, I\BPBI P.%
\BCBT {}\ \BBA {} Gray, W\BPBI G.%
\end{APACrefauthors}%
\unskip\
\newblock
\APACrefYearMonthDay{1984}{}{}.
\newblock
{\BBOQ}\APACrefatitle {One step integration methods with maximum stability
  regions} {One step integration methods with maximum stability
  regions}.{\BBCQ}
\newblock
\APACjournalVolNumPages{Math. Comput. Simulation}{26}{2}{87--92}.
\newblock
\begin{APACrefDOI} \doi{10.1016/0378-4754(84)90039-9} \end{APACrefDOI}
\PrintBackRefs{\CurrentBib}

\bibitem [\protect \citeauthoryear {%
Klein%
, Achatz%
, Bresch%
, Knio%
\BCBL {}\ \BBA {} Smolarkiewicz%
}{%
Klein%
\ \protect \BOthers {.}}{%
{\protect \APACyear {2010}}%
}]{%
klein_regime_2010}
\APACinsertmetastar {%
klein_regime_2010}%
\begin{APACrefauthors}%
Klein, R.%
, Achatz, U.%
, Bresch, D.%
, Knio, O\BPBI M.%
\BCBL {}\ \BBA {} Smolarkiewicz, P\BPBI K.%
\end{APACrefauthors}%
\unskip\
\newblock
\APACrefYearMonthDay{2010}{}{}.
\newblock
{\BBOQ}\APACrefatitle {Regime of validity of soundproof atmospheric flow
  models} {Regime of validity of soundproof atmospheric flow models}.{\BBCQ}
\newblock
\APACjournalVolNumPages{J. Atmos. Sci.}{67}{10}{3226--3237}.
\newblock
\begin{APACrefDOI} \doi{10.1175/2010JAS3490.1} \end{APACrefDOI}
\PrintBackRefs{\CurrentBib}

\bibitem [\protect \citeauthoryear {%
Laprise%
}{%
Laprise%
}{%
{\protect \APACyear {1992}}%
}]{%
laprise_euler_1992}
\APACinsertmetastar {%
laprise_euler_1992}%
\begin{APACrefauthors}%
Laprise, R.%
\end{APACrefauthors}%
\unskip\
\newblock
\APACrefYearMonthDay{1992}{}{}.
\newblock
{\BBOQ}\APACrefatitle {The {Euler} equations of motion with hydrostatic
  pressure as an independent variable} {The {Euler} equations of motion with
  hydrostatic pressure as an independent variable}.{\BBCQ}
\newblock
\APACjournalVolNumPages{Mon. Wea. Rev.}{120}{1}{197--207}.
\newblock
\begin{APACrefDOI} \doi{10.1175/1520-0493(1992)120<0197:TEEOMW>2.0.CO;2}
  \end{APACrefDOI}
\PrintBackRefs{\CurrentBib}

\bibitem [\protect \citeauthoryear {%
Lott%
, Woodward%
\BCBL {}\ \BBA {} Evans%
}{%
Lott%
\ \protect \BOthers {.}}{%
{\protect \APACyear {2015}}%
}]{%
lott_algorithmically_2015}
\APACinsertmetastar {%
lott_algorithmically_2015}%
\begin{APACrefauthors}%
Lott, P\BPBI A.%
, Woodward, C\BPBI S.%
\BCBL {}\ \BBA {} Evans, K\BPBI J.%
\end{APACrefauthors}%
\unskip\
\newblock
\APACrefYearMonthDay{2015}{}{}.
\newblock
{\BBOQ}\APACrefatitle {Algorithmically scalable block preconditioner for fully
  implicit shallow-water equations in {CAM}-{SE}} {Algorithmically scalable
  block preconditioner for fully implicit shallow-water equations in
  {CAM}-{SE}}.{\BBCQ}
\newblock
\APACjournalVolNumPages{Comp. Geosci.}{19}{1}{49--61}.
\newblock
\begin{APACrefDOI} \doi{10.1007/s10596-014-9447-6} \end{APACrefDOI}
\PrintBackRefs{\CurrentBib}

\bibitem [\protect \citeauthoryear {%
Marras%
\ \protect \BOthers {.}}{%
Marras%
\ \protect \BOthers {.}}{%
{\protect \APACyear {2016}}%
}]{%
marras_review_2016}
\APACinsertmetastar {%
marras_review_2016}%
\begin{APACrefauthors}%
Marras, S.%
, Kelly, J\BPBI F.%
, Moragues, M.%
, Müller, A.%
, Kopera, M\BPBI A.%
, Vázquez, M.%
\BDBL {}Jorba, O.%
\end{APACrefauthors}%
\unskip\
\newblock
\APACrefYearMonthDay{2016}{}{}.
\newblock
{\BBOQ}\APACrefatitle {A review of element-based {Galerkin} methods for
  numerical weather prediction: {Finite} elements, spectral elements, and
  discontinuous {Galerkin}} {A review of element-based {Galerkin} methods for
  numerical weather prediction: {Finite} elements, spectral elements, and
  discontinuous {Galerkin}}.{\BBCQ}
\newblock
\APACjournalVolNumPages{Arch. Comput. Methods Eng.}{23}{4}{673--722}.
\newblock
\begin{APACrefDOI} \doi{10.1007/s11831-015-9152-1} \end{APACrefDOI}
\PrintBackRefs{\CurrentBib}

\bibitem [\protect \citeauthoryear {%
Ogura%
\ \BBA {} Phillips%
}{%
Ogura%
\ \BBA {} Phillips%
}{%
{\protect \APACyear {1962}}%
}]{%
ogura_scale_1962}
\APACinsertmetastar {%
ogura_scale_1962}%
\begin{APACrefauthors}%
Ogura, Y.%
\BCBT {}\ \BBA {} Phillips, N\BPBI A.%
\end{APACrefauthors}%
\unskip\
\newblock
\APACrefYearMonthDay{1962}{}{}.
\newblock
{\BBOQ}\APACrefatitle {Scale analysis of deep and shallow convection in the
  atmosphere} {Scale analysis of deep and shallow convection in the
  atmosphere}.{\BBCQ}
\newblock
\APACjournalVolNumPages{J. Atmos. Sci.}{19}{2}{173--179}.
\newblock
\begin{APACrefDOI} \doi{10.1175/1520-0469(1962)019<0173:SAODAS>2.0.CO;2}
  \end{APACrefDOI}
\PrintBackRefs{\CurrentBib}

\bibitem [\protect \citeauthoryear {%
Rasch%
\ \protect \BOthers {.}}{%
Rasch%
\ \protect \BOthers {.}}{%
{\protect \APACyear {2018}}%
}]{%
rasch_overview_2018}
\APACinsertmetastar {%
rasch_overview_2018}%
\begin{APACrefauthors}%
Rasch, P\BPBI J.%
, Xie, S.%
, Ma, P\BHBI L.%
, Lin, W.%
, Wang, H.%
, Tang, Q.%
\BDBL {}Yang, Y.%
\end{APACrefauthors}%
\unskip\
\newblock
\APACrefYearMonthDay{2018}{}{}.
\newblock
{\BBOQ}\APACrefatitle {An overview of the atmospheric component of the {Energy}
  {Exascale} {Earth} {System} {Model}} {An overview of the atmospheric
  component of the {Energy} {Exascale} {Earth} {System} {Model}}.{\BBCQ}
\newblock
\APACjournalVolNumPages{J. Adv. Mod. Earth Sys.}{(under review)}{}{}.
\PrintBackRefs{\CurrentBib}

\bibitem [\protect \citeauthoryear {%
Rauscher%
\ \BBA {} Ringler%
}{%
Rauscher%
\ \BBA {} Ringler%
}{%
{\protect \APACyear {2014}}%
}]{%
rauscher_impact_2014}
\APACinsertmetastar {%
rauscher_impact_2014}%
\begin{APACrefauthors}%
Rauscher, S\BPBI A.%
\BCBT {}\ \BBA {} Ringler, T\BPBI D.%
\end{APACrefauthors}%
\unskip\
\newblock
\APACrefYearMonthDay{2014}{}{}.
\newblock
{\BBOQ}\APACrefatitle {Impact of variable-resolution meshes on midlatitude
  baroclinic eddies using {CAM}-{MPAS}-{A}} {Impact of variable-resolution
  meshes on midlatitude baroclinic eddies using {CAM}-{MPAS}-{A}}.{\BBCQ}
\newblock
\APACjournalVolNumPages{Mon. Wea. Rev.}{142}{11}{4256--4268}.
\newblock
\begin{APACrefDOI} \doi{10.1175/MWR-D-13-00366.1} \end{APACrefDOI}
\PrintBackRefs{\CurrentBib}

\bibitem [\protect \citeauthoryear {%
Reale%
, Achuthavarier%
, Fuentes%
, Putman%
\BCBL {}\ \BBA {} Partyka%
}{%
Reale%
\ \protect \BOthers {.}}{%
{\protect \APACyear {2017}}%
}]{%
reale_tropical_2017}
\APACinsertmetastar {%
reale_tropical_2017}%
\begin{APACrefauthors}%
Reale, O.%
, Achuthavarier, D.%
, Fuentes, M.%
, Putman, W\BPBI M.%
\BCBL {}\ \BBA {} Partyka, G.%
\end{APACrefauthors}%
\unskip\
\newblock
\APACrefYearMonthDay{2017}{}{}.
\newblock
{\BBOQ}\APACrefatitle {Tropical cyclones in the 7-km {NASA} global nature run
  for use in observing system simulation experiments} {Tropical cyclones in the
  7-km {NASA} global nature run for use in observing system simulation
  experiments}.{\BBCQ}
\newblock
\APACjournalVolNumPages{J. Atmospheric Ocean. Technol.}{34}{1}{73--100}.
\newblock
\begin{APACrefDOI} \doi{10.1175/JTECH-D-16-0094.1} \end{APACrefDOI}
\PrintBackRefs{\CurrentBib}

\bibitem [\protect \citeauthoryear {%
Rhoades%
\ \protect \BOthers {.}}{%
Rhoades%
\ \protect \BOthers {.}}{%
{\protect \APACyear {2018}}%
}]{%
rhoades_sensitivity_2018}
\APACinsertmetastar {%
rhoades_sensitivity_2018}%
\begin{APACrefauthors}%
Rhoades, A\BPBI M.%
, Ullrich, P\BPBI A.%
, Zarzycki, C\BPBI M.%
, Johansen, H.%
, Margulis, S\BPBI A.%
, Morrison, H.%
\BDBL {}Collins, W\BPBI D.%
\end{APACrefauthors}%
\unskip\
\newblock
\APACrefYearMonthDay{2018}{}{}.
\newblock
{\BBOQ}\APACrefatitle {Sensitivity of mountain hydroclimate simulations in
  variable-resolution {CESM} to microphysics and horizontal resolution}
  {Sensitivity of mountain hydroclimate simulations in variable-resolution
  {CESM} to microphysics and horizontal resolution}.{\BBCQ}
\newblock
\APACjournalVolNumPages{J. Adv. Mod. Earth Sys.}{10}{6}{1357--1380}.
\newblock
\begin{APACrefDOI} \doi{10.1029/2018MS001326} \end{APACrefDOI}
\PrintBackRefs{\CurrentBib}

\bibitem [\protect \citeauthoryear {%
Rokhzadi%
, Mohammadian%
\BCBL {}\ \BBA {} Charron%
}{%
Rokhzadi%
\ \protect \BOthers {.}}{%
{\protect \APACyear {2018}}%
}]{%
rokhzadi_optimally_2018}
\APACinsertmetastar {%
rokhzadi_optimally_2018}%
\begin{APACrefauthors}%
Rokhzadi, A.%
, Mohammadian, A.%
\BCBL {}\ \BBA {} Charron, M.%
\end{APACrefauthors}%
\unskip\
\newblock
\APACrefYearMonthDay{2018}{}{}.
\newblock
{\BBOQ}\APACrefatitle {An {Optimally} {Stable} and {Accurate} {Second}-{Order}
  {SSP} {Runge}-{Kutta} {IMEX} {Scheme} for {Atmospheric} {Applications}} {An
  {Optimally} {Stable} and {Accurate} {Second}-{Order} {SSP} {Runge}-{Kutta}
  {IMEX} {Scheme} for {Atmospheric} {Applications}}.{\BBCQ}
\newblock
\APACjournalVolNumPages{J. Adv. Mod. Earth Sys.}{10}{1}{18--42}.
\newblock
\begin{APACrefDOI} \doi{10.1002/2017MS001065} \end{APACrefDOI}
\PrintBackRefs{\CurrentBib}

\bibitem [\protect \citeauthoryear {%
Simmons%
\ \BBA {} Burridge%
}{%
Simmons%
\ \BBA {} Burridge%
}{%
{\protect \APACyear {1981}}%
}]{%
simmons_energy_1981}
\APACinsertmetastar {%
simmons_energy_1981}%
\begin{APACrefauthors}%
Simmons, A\BPBI J.%
\BCBT {}\ \BBA {} Burridge, D\BPBI M.%
\end{APACrefauthors}%
\unskip\
\newblock
\APACrefYearMonthDay{1981}{}{}.
\newblock
{\BBOQ}\APACrefatitle {An energy and angular-momentum conserving vertical
  finite-difference scheme and hybrid vertical coordinates} {An energy and
  angular-momentum conserving vertical finite-difference scheme and hybrid
  vertical coordinates}.{\BBCQ}
\newblock
\APACjournalVolNumPages{Mon. Wea. Rev.}{109}{4}{758--766}.
\newblock
\begin{APACrefDOI} \doi{10.1175/1520-0493(1981)109<0758:AEAAMC>2.0.CO;2}
  \end{APACrefDOI}
\PrintBackRefs{\CurrentBib}

\bibitem [\protect \citeauthoryear {%
Steyer%
, Vogl%
, Taylor%
\BCBL {}\ \BBA {} Guba%
}{%
Steyer%
\ \protect \BOthers {.}}{%
{\protect \APACyear {2019}}%
}]{%
steyer_efficient_2019}
\APACinsertmetastar {%
steyer_efficient_2019}%
\begin{APACrefauthors}%
Steyer, A.%
, Vogl, C\BPBI J.%
, Taylor, M.%
\BCBL {}\ \BBA {} Guba, O.%
\end{APACrefauthors}%
\unskip\
\newblock
\APACrefYearMonthDay{2019}{{\APACmonth{06}}}{}.
\newblock
{\BBOQ}\APACrefatitle {Efficient {IMEX} {Runge}-{Kutta} methods for
  nohydrostatic dynamics} {Efficient {IMEX} {Runge}-{Kutta} methods for
  nohydrostatic dynamics}.{\BBCQ}
\newblock
\APACjournalVolNumPages{SIAM J. Sci. Comput.}{(under review)}{}{preprint on
  arXiv: 1906.07219}.
\PrintBackRefs{\CurrentBib}

\bibitem [\protect \citeauthoryear {%
Taylor%
, Edwards%
, Thomas%
\BCBL {}\ \BBA {} Nair%
}{%
Taylor%
\ \protect \BOthers {.}}{%
{\protect \APACyear {2007}}%
}]{%
taylor_mass_2007}
\APACinsertmetastar {%
taylor_mass_2007}%
\begin{APACrefauthors}%
Taylor, M\BPBI A.%
, Edwards, J.%
, Thomas, S.%
\BCBL {}\ \BBA {} Nair, R.%
\end{APACrefauthors}%
\unskip\
\newblock
\APACrefYearMonthDay{2007}{}{}.
\newblock
{\BBOQ}\APACrefatitle {A mass and energy conserving spectral element
  atmospheric dynamical core on the cubed-sphere grid} {A mass and energy
  conserving spectral element atmospheric dynamical core on the cubed-sphere
  grid}.{\BBCQ}
\newblock
\APACjournalVolNumPages{J. Phys. Conf. Ser.}{78}{1}{012074}.
\newblock
\begin{APACrefDOI} \doi{10.1088/1742-6596/78/1/012074} \end{APACrefDOI}
\PrintBackRefs{\CurrentBib}

\bibitem [\protect \citeauthoryear {%
Taylor%
\ \BBA {} Fournier%
}{%
Taylor%
\ \BBA {} Fournier%
}{%
{\protect \APACyear {2010}}%
}]{%
taylor_compatible_2010}
\APACinsertmetastar {%
taylor_compatible_2010}%
\begin{APACrefauthors}%
Taylor, M\BPBI A.%
\BCBT {}\ \BBA {} Fournier, A.%
\end{APACrefauthors}%
\unskip\
\newblock
\APACrefYearMonthDay{2010}{}{}.
\newblock
{\BBOQ}\APACrefatitle {A compatible and conservative spectral element method on
  unstructured grids} {A compatible and conservative spectral element method on
  unstructured grids}.{\BBCQ}
\newblock
\APACjournalVolNumPages{J. Comp. Phys.}{229}{17}{5879--5895}.
\newblock
\begin{APACrefDOI} \doi{10.1016/j.jcp.2010.04.008} \end{APACrefDOI}
\PrintBackRefs{\CurrentBib}

\bibitem [\protect \citeauthoryear {%
{Taylor}%
\ \protect \BOthers {.}}{%
{Taylor}%
\ \protect \BOthers {.}}{%
{\protect \APACyear {2019}}%
}]{%
taylor_energy_2019}
\APACinsertmetastar {%
taylor_energy_2019}%
\begin{APACrefauthors}%
{Taylor}, M\BPBI A.%
, {Guba}, O.%
, {Steyer}, A.%
, {Ullrich}, P.%
, {Hall}, D.%
\BCBL {}\ \BBA {} {Eldred}, C.%
\end{APACrefauthors}%
\unskip\
\newblock
\APACrefYearMonthDay{2019}{{\APACmonth{08}}}{}.
\newblock
{\BBOQ}\APACrefatitle {An energy consistent discretization of the
  nonhydrostatic equations in primitive variables} {An energy consistent
  discretization of the nonhydrostatic equations in primitive
  variables}.{\BBCQ}
\newblock
\APACjournalVolNumPages{J. Adv. Mod. Earth Sys.}{(under review)}{}{preprint on
  arXiv: 1908.04430}.
\PrintBackRefs{\CurrentBib}

\bibitem [\protect \citeauthoryear {%
P.~Ullrich%
\ \BBA {} Jablonowski%
}{%
P.~Ullrich%
\ \BBA {} Jablonowski%
}{%
{\protect \APACyear {2012}}%
}]{%
ullrich_operator-split_2012}
\APACinsertmetastar {%
ullrich_operator-split_2012}%
\begin{APACrefauthors}%
Ullrich, P.%
\BCBT {}\ \BBA {} Jablonowski, C.%
\end{APACrefauthors}%
\unskip\
\newblock
\APACrefYearMonthDay{2012}{}{}.
\newblock
{\BBOQ}\APACrefatitle {Operator-split {Runge-–Kutta-–Rosenbrock} methods
  for nonhydrostatic atmospheric models} {Operator-split
  {Runge-–Kutta-–Rosenbrock} methods for nonhydrostatic atmospheric
  models}.{\BBCQ}
\newblock
\APACjournalVolNumPages{Mon. Wea. Rev.}{140}{4}{1257--1284}.
\newblock
\begin{APACrefDOI} \doi{10.1175/MWR-D-10-05073.1} \end{APACrefDOI}
\PrintBackRefs{\CurrentBib}

\bibitem [\protect \citeauthoryear {%
P\BPBI A.~Ullrich%
\ \protect \BOthers {.}}{%
P\BPBI A.~Ullrich%
\ \protect \BOthers {.}}{%
{\protect \APACyear {2017}}%
}]{%
ullrich_dcmip2016:_2017}
\APACinsertmetastar {%
ullrich_dcmip2016:_2017}%
\begin{APACrefauthors}%
Ullrich, P\BPBI A.%
, Jablonowski, C.%
, Kent, J.%
, Lauritzen, P\BPBI H.%
, Nair, R.%
, Reed, K\BPBI A.%
\BDBL {}Viner, K.%
\end{APACrefauthors}%
\unskip\
\newblock
\APACrefYearMonthDay{2017}{}{}.
\newblock
{\BBOQ}\APACrefatitle {{DCMIP}2016: a review of non-hydrostatic dynamical core
  design and intercomparison of participating models} {{DCMIP}2016: a review of
  non-hydrostatic dynamical core design and intercomparison of participating
  models}.{\BBCQ}
\newblock
\APACjournalVolNumPages{Geophys. Model Dev.}{10}{12}{4477--4509}.
\newblock
\begin{APACrefDOI} \doi{https://doi.org/10.5194/gmd-10-4477-2017}
  \end{APACrefDOI}
\PrintBackRefs{\CurrentBib}

\bibitem [\protect \citeauthoryear {%
P\BPBI A.~Ullrich%
\ \protect \BOthers {.}}{%
P\BPBI A.~Ullrich%
\ \protect \BOthers {.}}{%
{\protect \APACyear {2012}}%
}]{%
ullrich_dynamical_2012}
\APACinsertmetastar {%
ullrich_dynamical_2012}%
\begin{APACrefauthors}%
Ullrich, P\BPBI A.%
, Jablonowski, C.%
, Kent, J.%
, Lauritzen, P\BPBI H.%
, Nair, R\BPBI D.%
\BCBL {}\ \BBA {} Taylor, M\BPBI A.%
\end{APACrefauthors}%
\unskip\
\newblock
\APACrefYear{2012}.
\newblock
\APACrefbtitle {Dynamical core model intercomparison project ({DCMIP}) test
  case document} {Dynamical core model intercomparison project ({DCMIP}) test
  case document}.
\newblock
\APACrefnote{\url{https://earthsystemcog.org/site_media/docs/DCMIP-TestCaseDocument_v1.7.pdf}}
\PrintBackRefs{\CurrentBib}

\bibitem [\protect \citeauthoryear {%
P\BPBI A.~Ullrich%
, Reynolds%
, Guerra%
\BCBL {}\ \BBA {} Taylor%
}{%
P\BPBI A.~Ullrich%
\ \protect \BOthers {.}}{%
{\protect \APACyear {2018}}%
}]{%
ullrich_impact_2018}
\APACinsertmetastar {%
ullrich_impact_2018}%
\begin{APACrefauthors}%
Ullrich, P\BPBI A.%
, Reynolds, D\BPBI R.%
, Guerra, J\BPBI E.%
\BCBL {}\ \BBA {} Taylor, M\BPBI A.%
\end{APACrefauthors}%
\unskip\
\newblock
\APACrefYearMonthDay{2018}{}{}.
\newblock
{\BBOQ}\APACrefatitle {Impact and importance of hyperdiffusion on the spectral
  element method: {A} linear dispersion analysis} {Impact and importance of
  hyperdiffusion on the spectral element method: {A} linear dispersion
  analysis}.{\BBCQ}
\newblock
\APACjournalVolNumPages{J. Comp. Phys.}{375}{}{427--446}.
\newblock
\begin{APACrefDOI} \doi{10.1016/j.jcp.2018.06.035} \end{APACrefDOI}
\PrintBackRefs{\CurrentBib}

\bibitem [\protect \citeauthoryear {%
Vogl%
, Steyer%
, Reynolds%
, Ullrich%
\BCBL {}\ \BBA {} Woodward%
}{%
Vogl%
\ \protect \BOthers {.}}{%
{\protect \APACyear {2019}}%
{\protect \APACexlab {{\protect \BCnt {1}}}}}]{%
vogl_code_2019}
\APACinsertmetastar {%
vogl_code_2019}%
\begin{APACrefauthors}%
Vogl, C\BPBI J.%
, Steyer, A.%
, Reynolds, D\BPBI R.%
, Ullrich, P\BPBI A.%
\BCBL {}\ \BBA {} Woodward, C\BPBI S.%
\end{APACrefauthors}%
\unskip\
\newblock
\APACrefYearMonthDay{2019{\protect \BCnt {1}}}{}{}.
\newblock
\APACrefbtitle {Code {Archive}.} {Code {Archive}.}
\newblock
\begin{APACrefDOI} \doi{10.5281/zenodo.2559640} \end{APACrefDOI}
\PrintBackRefs{\CurrentBib}

\bibitem [\protect \citeauthoryear {%
Vogl%
, Steyer%
, Reynolds%
, Ullrich%
\BCBL {}\ \BBA {} Woodward%
}{%
Vogl%
\ \protect \BOthers {.}}{%
{\protect \APACyear {2019}}%
{\protect \APACexlab {{\protect \BCnt {2}}}}}]{%
vogl_stability_2019}
\APACinsertmetastar {%
vogl_stability_2019}%
\begin{APACrefauthors}%
Vogl, C\BPBI J.%
, Steyer, A.%
, Reynolds, D\BPBI R.%
, Ullrich, P\BPBI A.%
\BCBL {}\ \BBA {} Woodward, C\BPBI S.%
\end{APACrefauthors}%
\unskip\
\newblock
\APACrefYearMonthDay{2019{\protect \BCnt {2}}}{}{}.
\newblock
\APACrefbtitle {{Stability} {Plot} {Archive}.} {{Stability} {Plot} {Archive}.}
\newblock
\begin{APACrefDOI} \doi{10.5281/zenodo.3375044} \end{APACrefDOI}
\PrintBackRefs{\CurrentBib}

\bibitem [\protect \citeauthoryear {%
Wedi%
\ \protect \BOthers {.}}{%
Wedi%
\ \protect \BOthers {.}}{%
{\protect \APACyear {2015}}%
}]{%
wedi_modelling_2015}
\APACinsertmetastar {%
wedi_modelling_2015}%
\begin{APACrefauthors}%
Wedi, N.%
, Bauer, P.%
, Denoninck, W.%
, Diamantakis, M.%
, Hamrud, M.%
, Kuhnlein, C.%
\BDBL {}Smolarkiewicz, P.%
\end{APACrefauthors}%
\unskip\
\newblock
\APACrefYear{2015}.
\newblock
\APACrefbtitle {The modelling infrastructure of the {Integrated} {Forecasting}
  {System}: {Recent} advances and future challenges} {The modelling
  infrastructure of the {Integrated} {Forecasting} {System}: {Recent} advances
  and future challenges}.
\newblock
\APACrefnote{\url{https://www.ecmwf.int/sites/default/files/elibrary/2015/15259-modelling-infrastructure-integrated-forecasting-system-recent-advances-and-future-challenges.pdf}}
\PrintBackRefs{\CurrentBib}

\bibitem [\protect \citeauthoryear {%
Weller%
, Lock%
\BCBL {}\ \BBA {} Wood%
}{%
Weller%
\ \protect \BOthers {.}}{%
{\protect \APACyear {2013}}%
}]{%
weller_rungekutta_2013}
\APACinsertmetastar {%
weller_rungekutta_2013}%
\begin{APACrefauthors}%
Weller, H.%
, Lock, S\BHBI J.%
\BCBL {}\ \BBA {} Wood, N.%
\end{APACrefauthors}%
\unskip\
\newblock
\APACrefYearMonthDay{2013}{}{}.
\newblock
{\BBOQ}\APACrefatitle {{Runge-–Kutta} {IMEX} schemes for the horizontally
  explicit/vertically implicit ({HEVI}) solution of wave equations}
  {{Runge-–Kutta} {IMEX} schemes for the horizontally explicit/vertically
  implicit ({HEVI}) solution of wave equations}.{\BBCQ}
\newblock
\APACjournalVolNumPages{J. Comp. Phys.}{252}{}{365--381}.
\newblock
\begin{APACrefDOI} \doi{10.1016/j.jcp.2013.06.025} \end{APACrefDOI}
\PrintBackRefs{\CurrentBib}

\bibitem [\protect \citeauthoryear {%
Woodward%
\ \protect \BOthers {.}}{%
Woodward%
\ \protect \BOthers {.}}{%
{\protect \APACyear {2018}}%
}]{%
sundials-web}
\APACinsertmetastar {%
sundials-web}%
\begin{APACrefauthors}%
Woodward, C\BPBI S.%
, Reynolds, D\BPBI R.%
, Gardner, D\BPBI J.%
, Hindmarsh, A\BPBI C.%
, Balos, C\BPBI J.%
\BCBL {}\ \BBA {} Peles, S.%
\end{APACrefauthors}%
\unskip\
\newblock
\APACrefYearMonthDay{2018}{}{}.
\newblock
\APACrefbtitle {SUNDIALS Web Page.} {Sundials web page.}
\newblock
\APACrefnote{\url{https://computation.llnl.gov/projects/sundials/sundials-software}}
\PrintBackRefs{\CurrentBib}

\bibitem [\protect \citeauthoryear {%
Zarzycki%
\ \protect \BOthers {.}}{%
Zarzycki%
\ \protect \BOthers {.}}{%
{\protect \APACyear {2014}}%
}]{%
zarzycki_aquaplanet_2014}
\APACinsertmetastar {%
zarzycki_aquaplanet_2014}%
\begin{APACrefauthors}%
Zarzycki, C\BPBI M.%
, Levy, M\BPBI N.%
, Jablonowski, C.%
, Overfelt, J\BPBI R.%
, Taylor, M\BPBI A.%
\BCBL {}\ \BBA {} Ullrich, P\BPBI A.%
\end{APACrefauthors}%
\unskip\
\newblock
\APACrefYearMonthDay{2014}{}{}.
\newblock
{\BBOQ}\APACrefatitle {Aquaplanet experiments using {CAM}'s variable-resolution
  dynamical core} {Aquaplanet experiments using {CAM}'s variable-resolution
  dynamical core}.{\BBCQ}
\newblock
\APACjournalVolNumPages{J. Climate}{27}{14}{5481--5503}.
\newblock
\begin{APACrefDOI} \doi{10.1175/JCLI-D-14-00004.1} \end{APACrefDOI}
\PrintBackRefs{\CurrentBib}

\bibitem [\protect \citeauthoryear {%
Zerroukat%
, Wood%
\BCBL {}\ \BBA {} Staniforth%
}{%
Zerroukat%
\ \protect \BOthers {.}}{%
{\protect \APACyear {2006}}%
}]{%
zerroukat_parabolic_2006}
\APACinsertmetastar {%
zerroukat_parabolic_2006}%
\begin{APACrefauthors}%
Zerroukat, M.%
, Wood, N.%
\BCBL {}\ \BBA {} Staniforth, A.%
\end{APACrefauthors}%
\unskip\
\newblock
\APACrefYearMonthDay{2006}{}{}.
\newblock
{\BBOQ}\APACrefatitle {The Parabolic Spline Method (PSM) for conservative
  transport problems} {The parabolic spline method (psm) for conservative
  transport problems}.{\BBCQ}
\newblock
\APACjournalVolNumPages{Internat. J. Numer. Methods Fluids}{51}{11}{1297-1318}.
\newblock
\begin{APACrefDOI} \doi{10.1002/fld.1154} \end{APACrefDOI}
\PrintBackRefs{\CurrentBib}

\end{thebibliography}
%
% don't specify bibliographystyle
%%%%%%%%%%%%%%%%%%%%%%%%%%%%%%%%%%%%%%%%%%%%%%%

%Reference citation instructions and examples:
%
% Please use ONLY \cite and \citeA for reference citations.
% \cite for parenthetical references
% ...as shown in recent studies (Simpson et al., 2019)
% \citeA for in-text citations
% ...Simpson et al. (2019) have shown...
%
%
%...as shown by \citeA{jskilby}.
%...as shown by \citeA{lewin76}, \citeA{carson86}, \citeA{bartoldy02}, and \citeA{rinaldi03}.
%...has been shown \cite{jskilbye}.
%...has been shown \cite{lewin76,carson86,bartoldy02,rinaldi03}.
%... \cite <i.e.>[]{lewin76,carson86,bartoldy02,rinaldi03}.
%...has been shown by \cite <e.g.,>[and others]{lewin76}.
%
% apacite uses < > for prenotes and [ ] for postnotes
% DO NOT use other cite commands (e.g., \citet, \citep, \citeyear, \nocite, \citealp, etc.).
%

\end{document}